\documentclass[11pt]{amsart}

\usepackage{kbordermatrix}
\usepackage{graphicx}
\usepackage{wrapfig}
\usepackage{amssymb}
\usepackage{times}
\usepackage[top=20mm, bottom=20mm, left=20mm, right=20mm]{geometry}
\usepackage{color}
\usepackage{bm}
\usepackage{algorithmic}

\newcommand{\vct}[1]{\bm{\mathsf{#1}}}
\newcommand{\mtx}[1]{\bm{\mathsf{#1}}}

\newtheorem{thm}{Theorem}

\theoremstyle{definition}
\newtheorem{remark}{Remark}

\newcommand{\lsp}{\vspace{3mm}}

\newcommand{\vtwo}[2]{\left[\begin{array}{c} #1 \\ #2 \end{array}\right]}

\newcommand{\mtwo}[4]{\left[\begin{array}{cc}    #1 & #2 \\ #3 & #4  \end{array}\right]}

\newcommand{\randUTV}{randUTV}

\begin{document}

\begin{center}
\textbf{\large{randUTV: A blocked randomized algorithm for computing a rank-revealing UTV factorization}}

\vspace{5mm}

\textit{
P.G.~Martinsson\footnotemark[1],
G.~Quintana-Ort\'{\i}\footnotemark[2],
N.~Heavner\footnotemark[1]}

\vspace{5mm}

\begin{minipage}{140mm}
\textbf{Abstract:}
This manuscript describes the randomized algorithm \randUTV{}
for computing a so called UTV factorization efficiently.
Given a matrix $\mtx{A}$,
the algorithm computes a factorization $\mtx{A} = \mtx{U}\mtx{T}\mtx{V}^{*}$,
where $\mtx{U}$ and $\mtx{V}$ have orthonormal columns,
and $\mtx{T}$ is triangular (either upper or lower, whichever is preferred).
The algorithm \randUTV{} is developed primarily
to be a fast and easily parallelized alternative
to algorithms for computing the Singular Value Decomposition (SVD).
\randUTV{} provides accuracy very close to that of the SVD for problems
such as low-rank approximation, solving ill-conditioned linear systems,
determining bases for various subspaces associated with the matrix, etc.
Moreover, \randUTV{} also produces
highly accurate approximations to the singular values of $\mtx{A}$.
Unlike the SVD,
the randomized algorithm proposed builds a UTV factorization in an
incremental, single-stage, and non-iterative way, making it
possible to halt the factorization process once a specified tolerance has been met.
Numerical experiments comparing the accuracy and speed of \randUTV{}
to the SVD are presented.
These experiments demonstrate that in comparison to column pivoted QR, which is another
factorization that is often used as a relatively economic alternative to the SVD,
\randUTV{} compares favorably in terms of speed while providing far higher
accuracy.
\end{minipage}

\end{center}

\footnotetext[1]{Department of Applied Mathematics, University of Colorado at Boulder, 526 UCB, Boulder, CO 80309-0526, USA}

\footnotetext[2]{Depto.~de Ingenier\'{\i}a y Ciencia de Computadores, Universitat Jaume I, 12.071 Castell\'on, Spain}

\lsp

\section{Introduction}

\subsection{Overview}

Given an $m\times n$ matrix $\mtx{A}$, the so called ``UTV decomposition'' \cite[p.~400]{1998_stewart_volume1} takes the form
\begin{equation}
\label{eq:defUTVpre}
\begin{array}{ccccccccccc}
\mtx{A} &=& \mtx{U} & \mtx{T} & \mtx{V}^{*},\\
m\times n && m\times m & m\times n & n\times n
\end{array}
\end{equation}
where $\mtx{U}$ and $\mtx{V}$ are unitary matrices, and
$\mtx{T}$ is a triangular matrix (either lower or upper triangular).
The UTV decomposition can be viewed as a generalization of other standard factorizations
such as, e.g., the \textit{Singular Value Decomposition (SVD)} or the
\textit{Column Pivoted QR decomposition (CPQR).} (To be precise, the SVD is the special case
where $\mtx{T}$ is diagonal, and the CPQR is the special case where $\mtx{V}$ is a permutation matrix.)
The additional flexibility inherent in the UTV
decomposition enables the design of efficient updating procedures,
see \cite[Ch.~5, Sec.~4]{1998_stewart_volume1} and \cite{1999_hansen_UTVtools,1994_stewart_UTV,1995_elden_downdate_UTV}.
In this manuscript, we describe a randomized algorithm we call \randUTV{} that
exploits the additional flexibility provided by the UTV format to build a factorization
algorithm that combines some of the most desirable properties of standard algorithms for
computing the SVD and CPQR factorizations.

Specific advantages of the proposed algorithm include:
(i) \randUTV{} provides close to optimal low-rank approximation in the sense that
for $k = 1,2,\dots,\min(m,n)$, the output factors $\mtx{U}$, $\mtx{V}$, and $\mtx{T}$ satisfy
$$
\|\mtx{A} - \mtx{U}(:,1:k)\mtx{T}(1:k,:)\mtx{V}^{*}\| \approx \inf\{\|\mtx{A} - \mtx{B}\|\,\colon\,\mtx{B}\mbox{ has rank }k\}.
$$
In particular, \randUTV{} is much better at low-rank approximation than CPQR.
A related advantage of \randUTV{} is that the diagonal elements of $\mtx{T}$
provide excellent approximations to the singular values of $\mtx{A}$, cf.~Figure \ref{fig:diags}.
(ii) The algorithm \randUTV{} builds the factorization (\ref{eq:defUTVpre}) incrementally,
which means that when it is applied to a matrix of numerical rank $k$, the algorithm can be
stopped early and incur an overall cost of $O(mnk)$. Observe that standard algorithms for computing
the SVD do not share this property. The CPQR algorithm can be stopped in this fashion, but
typically leads to substantially suboptimal low-rank approximation.
(iii) The algorithm \randUTV{} is \textit{blocked}. For a given block size $b$,
\randUTV{} processes $b$ columns of $\mtx{A}$ at a time. Moreover, the vast majority
of flops expended by \randUTV{} are used in matrix-matrix multiplications involving $\mtx{A}$
and thin matrices with $b$ columns. This leads to high computational speeds, in particular when
the algorithm is executed on multicore CPUs, GPUs, and
distributed-memory architectures.
(iv) The algorithm \randUTV{} is not an iterative algorithm. In this regard, it is closer
to the CPQR than to standard SVD algorithms, which substantially simplifies software optimization.

\begin{wrapfigure}{r}{100mm}
\centering
\includegraphics[width=80mm]{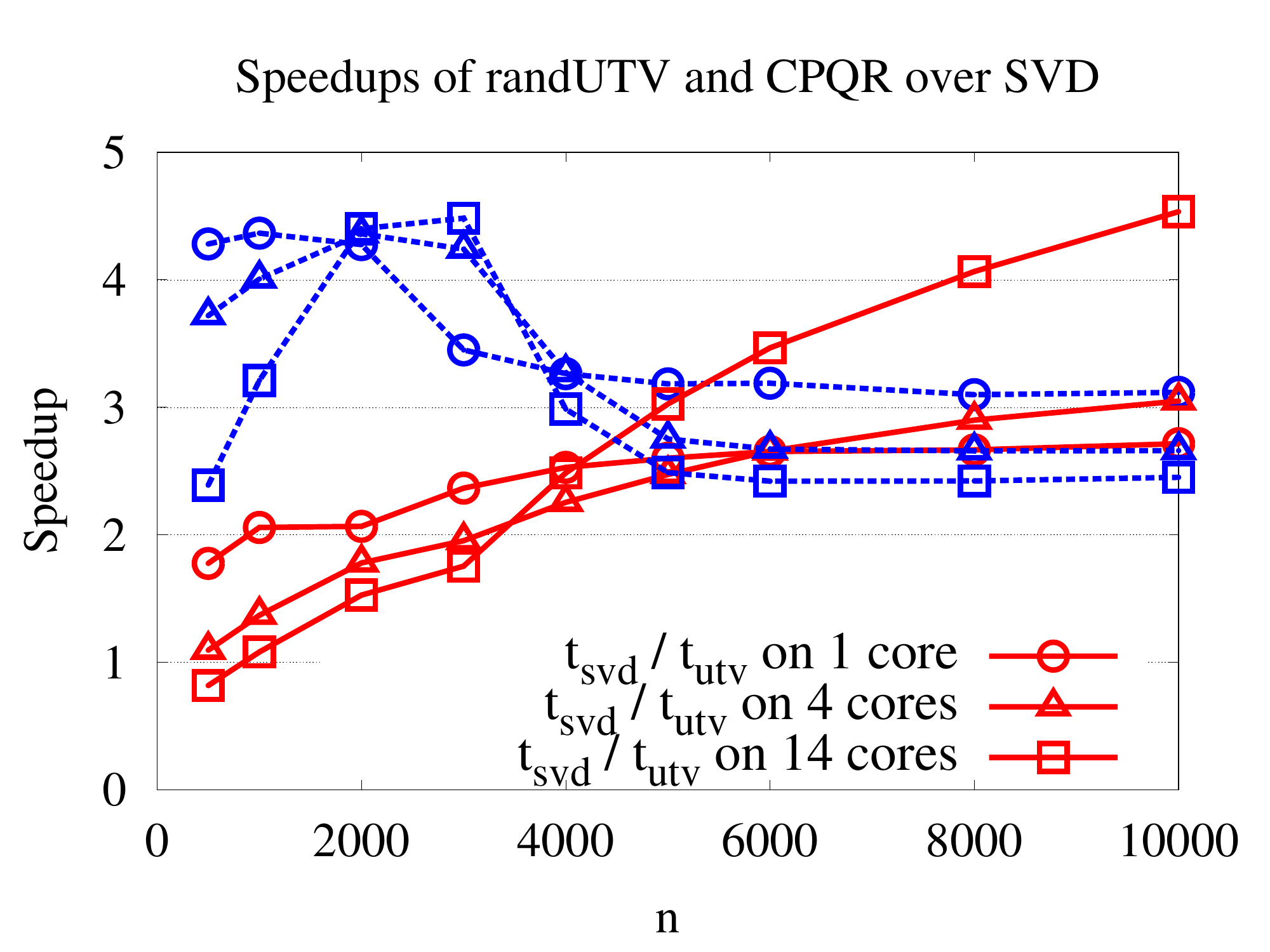}
\caption{Solid lines show the acceleration of \randUTV{} with the parameter choice
$q=1$ (cf.~Section \ref{sec:right}) compared to Intel's MKL SVD (\texttt{dgesvd}).
Dashed lines show the analogous acceleration of CPQR over SVD.}
\label{fig:speedup}
\end{wrapfigure}

Section \ref{sec:num} describes the results from several numerical experiments that illustrate
the accuracy and speed of \randUTV{} in comparison to standard techniques for computing
the SVD and CPQR factorizations. As a preview, we show in Figure \ref{fig:speedup} that
\randUTV{} executes far faster than a highly optimized implementation of the LAPACK
function \texttt{dgesvd} for computing the SVD (we compare to the Intel MKL implementation).
Figure \ref{fig:speedup} also includes lines that indicate how much faster CPQR is over the
SVD, since the CPQR is often proposed as an economical alternative to the SVD for low rank
approximation.
Although the comparison between CPQR and \randUTV{} is slightly unfair
since \randUTV{} is far more accurate than CPQR,
we see that CPQR is faster for small matrices, but
\randUTV{} becomes competitive as the matrix size increases.
More importantly,
the relative speed of \randUTV{} increases greatly as the number of processors increases.
In other words, in modern computing environments, \randUTV{} is \textit{both} faster \textit{and}
far better at revealing the numerical rank of $\mtx{A}$ than column pivoted QR.

\subsection{A randomized algorithm for computing the UTV decomposition}

The algorithm we propose is blocked for computational efficiency.
For concreteness, let us assume that $m\geq n$, and that an \textit{upper}
triangular factor $\mtx{T}$ is sought.
With $b$ a block size,
\randUTV{} proceeds through approximately $n/b$ steps,
where at the $i$'th step the $i$'th block of $b$ columns of
$\mtx{A}$ is driven to upper triangular form, as illustrated in Figure
\ref{fig:block_cartoon}.

In the first iteration, \randUTV{}
uses a randomized subspace iteration inspired by
\cite{2009_szlam_power,2011_martinsson_randomsurvey} to build
two sets of $b$ orthonormal vectors that approximately span the spaces spanned
by the $b$ dominant left and right singular vectors of $\mtx{A}$, respectively.
These basis vectors form the first $b$ columns of two unitary ``transformation
matrices'' $\mtx{U}^{(1)}$ and $\mtx{V}^{(1)}$. We use these to build a new matrix
$$
\mtx{A}^{(1)} = \bigl(\mtx{U}^{(1)}\bigr)^{*}\,\mtx{A}\,\mtx{V}^{(1)}
$$
that has a blocked structure as follows:
$$
\mtx{A}^{(1)} =
\left[\begin{array}{cc}
\mtx{A}^{(1)}_{11} & \mtx{A}^{(1)}_{12} \\
\mtx{0}            & \mtx{A}^{(1)}_{22}
\end{array}\right] =
\raisebox{-9mm}{\includegraphics[height=18mm]{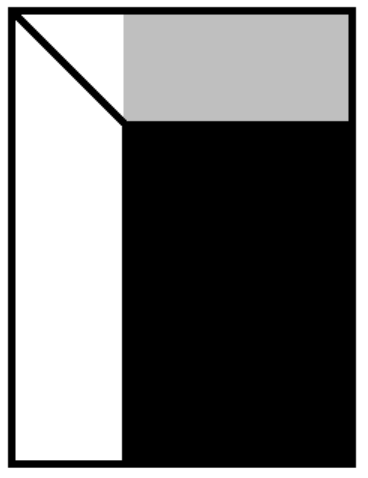}}.
$$
In other words, the top left $b\times b$ block is diagonal, and the bottom left
block is zero. Ideally, we would want the top right block to also be zero, and
this is what we would obtain if we used the exact left and right singular vectors
in the transformation.
The randomized sampling does not exactly identity the
correct subspaces, but it does to high enough accuracy that the elements in the
top right block have very small moduli. For purposes of low-rank approximation,
such a format strikes an attractive balance between computational efficiency
and close to optimal accuracy. We will demonstrate that
$\|\mtx{A}_{22}^{(1)}\| \approx \inf\{\|\mtx{A} - \mtx{B}\|\,\colon\,\mtx{B}\mbox{ has rank }b\}$,
and that the diagonal entries of $\mtx{A}^{(1)}_{11}$ form accurate approximations to the first $b$
singular values of $\mtx{A}$.

Once the first step has been completed, the second step applies the same procedure
to the remaining block $\mtx{A}^{(1)}_{22}$, which has size $(m-b) \times (n-b)$,
and then continues in the same fashion to process all remaining blocks.

\subsection{Relationship to earlier work}

The UTV factorization was introduced and popularized by G.W.~Stewart in a sequence of papers,
including \cite{1992_stewart_subspace_tracking,1993_stewart_update_ULV,1994_stewart_UTV,1999_stewart_QLP},
and the text book chapter \cite[p.~400]{1998_stewart_volume1}. Among these,
\cite{1999_stewart_QLP} is of particular relevance as it discusses explicitly how
the UTV decomposition can be used for low-rank approximation and for finding approximations
to the singular values of a matrix; in Section \ref{sec:errors} we compare the accuracy of
the method in \cite{1999_stewart_QLP} to the accuracy of \randUTV{}.
Of relevance here is also \cite{1993_stewart_blockQR_SVD},
which describes deterministic iterative methods
for driving a triangular matrix towards diagonality.
Another path towards better rank revelation and more accurate rank estimation is described in
\cite{1997_pchansen_low_rank_UTV}.
A key advantage of the UTV factorization is the ease with which it can be updated, as discussed in,
e.g., \cite{1996_barlow_downdating_ULV,2003_barlow_UTV,1995_elden_downdate_UTV}. Implementation aspects
are discussed in \cite{1999_hansen_UTVtools}.

The work presented here relies crucially on previous work on randomized subspace iteration for
approximating the linear spaces spanned by groups of singular vectors, including
\cite{2009_szlam_power,2011_martinsson_randomsurvey,2011_martinsson_random1,2006_martinsson_random1_orig,2011_halko_largedata}.
This prior body of literature focused on the task of computing \textit{partial} factorizations
of matrices. More recently, it was observed \cite{2015_blockQR,2015_blockQR_inreview,2015_blockQR_ming} that
analogous techniques can be utilized to accelerate methods for computing \textit{full} factorizations
such as the CPQR. The gain in speed is attained from \textit{blocking} of the algorithms, rather than by reducing
the asymptotic flop count. An alternative approach to using randomization was described in
\cite{2007_demmel_fast_linear_algebra_is_stable}.

\subsection{Outline of paper}

Section \ref{sec:prel} introduces notation, and lists some relevant existing results that we need.
Section \ref{sec:overview} provides a high-level description of the proposed algorithm.
Section \ref{sec:singlestep} describes in detail how to drive one block of columns to upper
triangular form, which forms one step of the blocked algorithm.
Section \ref{sec:randUTV} describes the whole multistep algorithm, discusses some implementation
issues, and provides an estimate of the asymptotic flop count.
Section \ref{sec:num} gives the results of several numerical experiments that illustrate the speed
and accuracy of the proposed method.
Section \ref{sec:code} describes availability of software.

\section{Preliminaries}
\label{sec:prel}

This section introduces our notation, and reviews some established techniques that
will be needed. The material in Sections \ref{sec:notation}--\ref{sec:tallthin}
is described in any standard text on numerical linear algebra
(e.g.~\cite{golub,1998_stewart_volume1,trefethen1997numerical}). The material in
Section \ref{sec:randrange} on randomized algorithms is described in further detail
in the survey \cite{2011_martinsson_randomsurvey} and the lecture notes
\cite{2016_martinsson_randomized_notes}.

\subsection{Basic notation}
\label{sec:notation}

Throughout this manuscript, we measure vectors in $\mathbb{R}^{n}$ using their
Euclidean norm. The default norm for matrices will be the corresponding
operator norm $\|\mtx{A}\| = \sup\{\|\mtx{A}\vct{x}\|\,\colon\,\|\vct{x}\|=1\}$,
although we will sometimes also use the Frobenius norm
$\|\mtx{A}\|_{\rm Fro} = \bigl(\sum_{i,j}|\mtx{A}(i,j)|^{2}\bigr)^{1/2}$.
We use the notation of Golub and Van Loan \cite{golub} to specify
submatrices: If $\mtx{B}$ is an $m\times n$ matrix,
and $I = [i_{1},\,i_{2},\,\dots,\,i_{k}]$ and $J = [j_{1},\,j_{2},\,\dots,\,j_{\ell}]$
are index vectors, then $\mtx{B}(I,J)$ denotes the corresponding $k\times \ell$
submatrix..
We let $\mtx{B}(I,:)$ denote the matrix $\mtx{B}(I,[1,\,2,\,\dots,\,n])$,
and define $\mtx{B}(:,J)$ analogously.
$\mtx{I}_{n}$ denotes the $n\times n$ identity matrix, and $\mtx{0}_{m,n}$ is the
$m\times n$ zero matrix.
The transpose of $\mtx{B}$ is denoted $\mtx{B}^{*}$, and we say that an $m\times n$ matrix $\mtx{U}$
is \textit{orthonormal} if its columns are orthonormal, so that $\mtx{U}^{*}\mtx{U} = \mtx{I}_{n}$.
A square orthonormal matrix is said to be \textit{unitary.}

\subsection{The Singular Value Decomposition (SVD)}

Let $\mtx{A}$ be an $m\times n$ matrix and set $r = \min(m,n)$.
Then the SVD of $\mtx{A}$ takes the form
\begin{equation}
\label{eq:defSVD}
\begin{array}{ccccccccccc}
\mtx{A} &=& \mtx{U} & \mtx{D} & \mtx{V}^{*},\\
m\times n && m\times r & r\times r & r\times n
\end{array}
\end{equation}
where matrices $\mtx{U}$ and $\mtx{V}$ are orthonormal,
and $\mtx{D}$ is diagonal. We have
$\mtx{U} = \bigl[\vct{u}_{1}\ \vct{u}_{2}\ \dots\ \vct{u}_{r}\bigr]$,
$\mtx{V} = \bigl[\vct{v}_{1}\ \vct{v}_{2}\ \dots\ \vct{v}_{r}\bigr]$,
and
$\mtx{D} = \mbox{diag}(\sigma_{1},\,\sigma_{2},\,\dots,\,\sigma_{r})$,
where $\{\vct{u}_{j}\}_{j=1}^{r}$ and $\{\vct{v}_{j}\}_{j=1}^{r}$ are the
left and right singular vectors of $\mtx{A}$, respectively, and
$\{\sigma_{j}\}_{j=1}^{r}$ are the singular values of $\mtx{A}$.
The singular values are ordered so that
$\sigma_{1} \geq \sigma_{2} \geq \dots \geq \sigma_{r} \geq 0$.

A principal advantage of the SVD is that it furnishes an explicit solution to
the low rank approximation problem. To be precise, let us for $k = 1,2,\dots,r$
define the rank-$k$ matrix
$$
\mtx{A}_{k} = \mtx{U}(:,1:k)\,\mtx{D}(1:k,1:k)\,\mtx{V}(:,1:k)^{*} =
\sum_{j=1}^{k}\sigma_{j}\,\vct{u}_{j}\,\vct{v}_{j}^{*}.
$$
Then, the Eckart-Young theorem asserts that
$$
\|\mtx{A} - \mtx{A}_{k}\| = \inf\{\|\mtx{A} - \mtx{B}\|\,\colon\,\mtx{B}\mbox{ has rank }k\}.
$$
A disadvantage of the SVD is that the cost to compute the full SVD is $O(n^{3})$ for a square
matrix, and $O(mnr)$ for a rectangular matrix, with large pre-factors. Moreover, standard techniques
for computing the SVD are challenging to parallelize, and cannot readily be modified to compute
partial factorizations.

\subsection{The Column Pivoted QR (CPQR) decomposition}
\label{sec:CPQR}

Let $\mtx{A}$ be an $m\times n$ matrix and set $r = \min(m,n)$.
Then, the CPQR decomposition of $\mtx{A}$ takes the form
\begin{equation}
\label{eq:defCPQR}
\begin{array}{ccccccccccc}
\mtx{A} &=& \mtx{Q} & \mtx{R} & \mtx{P}^{*},\\
m\times n && m\times r & r\times n & n\times n
\end{array}
\end{equation}
where $\mtx{Q}$ is orthonormal, $\mtx{R}$ is upper triangular, and
$\mtx{P}$ is a permutation matrix. The permutation matrix $\mtx{P}$ is typically
chosen to ensure monotonic decay in magnitude of the diagonal entries of $\mtx{R}$
so that $|\mtx{R}(1,1)| \geq |\mtx{R}(2,2)| \geq |\mtx{R}(3,3)| \geq \cdots$. The
factorization (\ref{eq:defCPQR}) is commonly computed using the \textit{Householder
QR algorithm} \cite[Sec.~5.2]{golub}, which is exceptionally stable.

An advantage of the CPQR is that it is computed via an incremental algorithm that
can be halted to produce a partial factorization once any given tolerance has been
met. A disadvantage is that it is not ideal for low-rank approximation.
In the typical case, the error incurred is noticeably worse than what you get from
the SVD but not disastrously so. However, there exist matrices for which CPQR leads
to very suboptimal approximation errors \cite{1966_kahan_NLA}.
Specialized pivoting strategies have been developed that can in some circumstances
improve the low-rank approximation property, resulting in so called ``rank-revealing
QR factorizations'' \cite{1987_chan_RRQR,gu1996}.

The classical column pivoted Householder QR factorization algorithm drives the matrix
$\mtx{A}$ to upper triangular form via a sequence of $r-1$ rank-one updates
and column pivotings.
Due to the column pivoting performed after each update,
it is difficult to \textit{block}, making it hard to achieve high computational
efficiency on modern processors \cite{DDSV}. It has recently been
demonstrated that randomized methods can be used to resolve this difficulty \cite{2015_blockQR,2015_blockQR_inreview,2015_blockQR_ming}.
However, we do not yet have rigorous theory backing up such randomized techniques, and they
have not yet been incorporated into standard software packages.

\subsection{Efficient factorizations of tall and thin matrices}
\label{sec:tallthin}

The algorithms proposed in this manuscript rely crucially on the fact
that factorizations of ``tall and thin'' matrices can be computed efficiently. To
be precise, suppose that we are given a matrix $\mtx{B}$ of size $m\times b$, where $m \gg b$,
and that we seek to compute an \textit{unpivoted} QR factorization
\begin{equation}
\label{eq:unpivQR}
\begin{array}{ccccccccccccccccc}
\mtx{B} &=& \mathbb{Q} & \mtx{R}, \\
m\times b && m\times m & m\times b
\end{array}
\end{equation}
where $\mathbb{Q}$ is unitary, and $\mtx{R}$ is upper triangular.
When the Householder
QR factorization procedure is used to compute the factorization (\ref{eq:unpivQR}), the
matrix $\mathbb{Q}$ is formed as a product of $b$ so called ``Householder reflectors''
\cite[Sec.~5.2]{golub}, and can be written in the form
\begin{equation}
\label{eq:WYrepresentation}
\begin{array}{cccccccccccccc}
\mathbb{Q} &=& \mtx{I} &+& \mtx{W} & \mtx{Y},\\
m\times m && m\times m && m\times b & b\times m
\end{array}
\end{equation}
for some matrices $\mtx{W}$ and $\mtx{Y}$ that can be readily computed given the
$b$ Householder vectors that are formed in the QR factorization procedure
\cite{BiVL87,ScVL89,Joffrain:2006:AHT:1141885.1141886}.
As a consequence, we need only $O(mb)$ words of storage to store $\mathbb{Q}$, and
we can apply $\mathbb{Q}$ to an $m\times n$ matrix using $\sim 2mnb$ flops.
In this manuscript,
the different typeface in $\mathbb{Q}$ is used as a reminder that this is a matrix
that can be stored and applied efficiently.

Next, suppose that we seek to compute the SVD of the tall and thin matrix $\mtx{B}$. This
can be done efficiently in a two-step procedure, where the first step is to compute the
unpivoted QR factorization (\ref{eq:unpivQR}). Next, let $\mtx{R}_{\rm small}$ denote
the top $b\times b$ block of $\mtx{R}$ so that
\begin{equation}
\label{eq:dallas}
\mtx{R} = \vtwo{\mtx{R}_{\rm small}}{\mtx{0}_{m-b,b}}.
\end{equation}
Then, compute the SVD of $\mtx{R}_{\rm small}$ to obtain the factorization
\begin{equation}
\label{eq:smallSVD}
\begin{array}{cccccccccccccccccc}
\mtx{R}_{\rm small} &=& \mtx{U}_{\rm small} & \mtx{D}_{\rm small} & \mtx{V}_{\rm small}^{*}.\\
b \times b && b \times b & b \times b & b \times b
\end{array}
\end{equation}
Combining (\ref{eq:unpivQR}), (\ref{eq:dallas}), and (\ref{eq:smallSVD}), we obtain the factorization
\begin{equation}
\label{eq:tallSVD}
\begin{array}{ccccccccccccccccc}
\mtx{B} &=& \mathbb{Q} & \mtwo{\mtx{U}_{\rm small}}{\mtx{0}}{\mtx{0}}{\mtx{I}_{m-b}} &
\vtwo{\mtx{D}_{\rm small}}{\mtx{0}_{m-b,b}} & \mtx{V}_{\rm small}^{*}, \\
m\times b && m\times m & m\times m & m\times b &b\times b
\end{array}
\end{equation}
which we recognize as a singular value decomposition of $\mtx{B}$. Observe that the cost of
computing this factorization is $O(mb^{2})$, and that only $O(mb)$ words of storage are required,
despite the fact that the matrix of left singular vectors is ostensibly an $m\times m$ dense matrix.

\begin{remark}
We mentioned in Section \ref{sec:CPQR} that it is challenging to achieve high
performance when implementing column pivoted QR on modern processors. In contrast,
\textit{unpivoted} QR is highly efficient, since this algorithm can readily be blocked
(see, e.g., Figures 4.1 and 4.2 in \cite{2015_blockQR_inreview}).
\end{remark}

\subsection{Randomized power iterations}
\label{sec:randrange}

This section summarizes key results of \cite{2009_szlam_power,2011_martinsson_randomsurvey,2016_martinsson_randomized_notes}
on randomized algorithms for constructing sets of orthonormal vectors that approximately
span the row or column spaces of a given matrix. To be precise,
let $\mtx{A}$ be an $m\times n$ matrix, let $b$ be an integer such that $1 \leq b < \min(m,n)$,
and suppose that we seek to find an $n\times b$ orthonormal matrix $\mtx{Q}$ such that
$$
\|\mtx{A} - \mtx{A}\mtx{Q}\mtx{Q}^{*}\| \approx \inf\{\|\mtx{A} - \mtx{B}\|\,\colon\,\mtx{B}\mbox{ has rank }b\}.
$$
Informally, the columns of $\mtx{Q}$ should approximately span the same space as
the dominant $b$ right singular vectors of $\mtx{A}$. This is a task that is very
well suited for subspace iteration (see \cite[Sec.~4.4.3]{1997_demmel_ANLA} and \cite{bathe1971solution}),
in particular when the starting matrix is a Gaussian random matrix
\cite{2009_szlam_power,2011_martinsson_randomsurvey}. With $q$ a (small) integer
denoting the number of steps taken in the power iteration, the following algorithm
leads to close to optimal results:
\begin{enumerate}
\item Draw a Gaussian random matrix $\mtx{G}$ of size $m\times b$.
\item Form a ``sampling matrix'' $\mtx{Y}$ of size $n\times b$ via $\mtx{Y} = \bigl(\mtx{A}^{*}\mtx{A}\bigr)^{q}\mtx{A}^{*}\mtx{G}$.
\item Orthonormalize the columns of $\mtx{Y}$ to form the matrix $\mtx{Q}$.
\end{enumerate}
Observe that in step (2), the matrix $\mtx{Y}$ is computed via alternating application
of $\mtx{A}^{*}$ and $\mtx{A}$ to a tall thin matrix with $b$ columns. In some situations,
orthonormalization is required between applications to avoid loss of accuracy due to
floating point arithmetic.
In \cite{2009_szlam_power,2011_martinsson_randomsurvey} it is demonstrated that by
using a Gaussian random matrix as the starting point, it is often sufficient to take
just a few steps of power iteration, say $q=1$ or $2$, or even $q=0$.

\begin{remark}[Over-sampling]
\label{remark:oversampling}
In the analysis of power iteration with a Gaussian random matrix as the starting point,
it is common to draw a few ``extra'' samples. In other words, one picks a small integer
$p$ representing the amount of over-sampling done, say $p=5$ or $p=10$, and starts with
a Gaussian matrix of size $m \times (b+p)$.
This results in an orthonormal (ON) matrix $\mtx{Q}$ of
size $n\times (b+p)$.
Then, with probability almost 1, the
error $\|\mtx{A} - \mtx{A}\mtx{Q}\mtx{Q}^{*}\|$ is close to the minimal error
in rank-$b$ approximation in both spectral and Frobenius norm \cite[Sec.~10]{2011_martinsson_randomsurvey}.
When no over-sampling is used, one risks losing some accuracy in the last couple of modes of the singular
value decomposition. However, our experience shows that in the context of the present
article, this loss is hardly noticeable.
\end{remark}

\begin{remark}[RSVD]
\label{remark:RSVD}
The randomized range finder described in this section is simply the
first stage in the two-stage ``Randomized SVD (RSVD)'' procedure for
computing an approximate
rank-$b$ SVD of a given matrix $\mtx{A}$ of size $m\times n$. To wit, suppose
that we have used the randomized range finder to build an ON matrix $\mtx{Q}$
of size $n\times b$ such that $\mtx{A} = \mtx{A}\mtx{Q}\mtx{Q}^{*} + \mtx{E}$,
for some some error matrix $\mtx{E}$.
Then, we can compute an approximate factorization
\begin{equation}
\label{eq:panda1}
\begin{array}{ccccccccccc}
\mtx{A} &=& \mtx{U} & \mtx{D} & \mtx{V}^{*} &+& \mtx{E},\\
m\times n && m\times b & b\times b & b \times n && m\times n
\end{array}
\end{equation}
where $\mtx{U}$ and $\mtx{V}$ are orthonormal, and $\mtx{D}$ is diagonal, via the
following steps (which together form the ``second stage'' of RSVD): (1) Set $\mtx{B} = \mtx{A}\mtx{Q}$ so that $\mtx{A}\mtx{Q}\mtx{Q}^{*} = \mtx{B}\mtx{Q}^{*}$. (2) Compute a full
SVD of the small matrix $\mtx{B}$ so that
\begin{equation}
\label{eq:panda2}
\begin{array}{ccccccccccc}
\mtx{B} &=& \mtx{U} & \mtx{D} & \hat{\mtx{V}}^{*}.\\
m\times b && m\times b & b\times b & b\times b
\end{array}
\end{equation}
(3) Set $\mtx{V} = \mtx{Q}\hat{\mtx{V}}$. Observe that these last three steps
are exact up to floating point arithmetic, so the error in (\ref{eq:panda1}) is
exactly the same as the error in the range finder alone:
$\mtx{E} = \mtx{A} - \mtx{A}\mtx{Q}\mtx{Q}^{*} = \mtx{A} - \mtx{U}\mtx{D}\mtx{V}^{*}$.
\end{remark}

\section{An overview of the randomized UTV factorization algorithm}
\label{sec:overview}

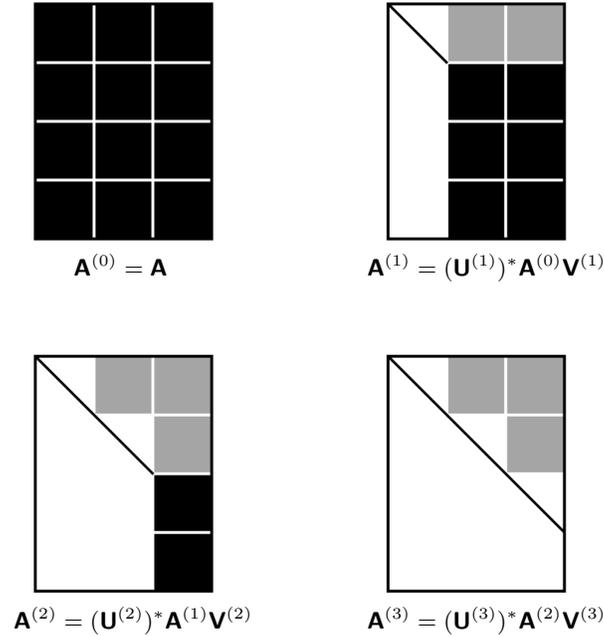
\begin{figure}
\setlength{\unitlength}{1mm}
\begin{picture}(78,85)
\put(000,05){\includegraphics[height=80mm]{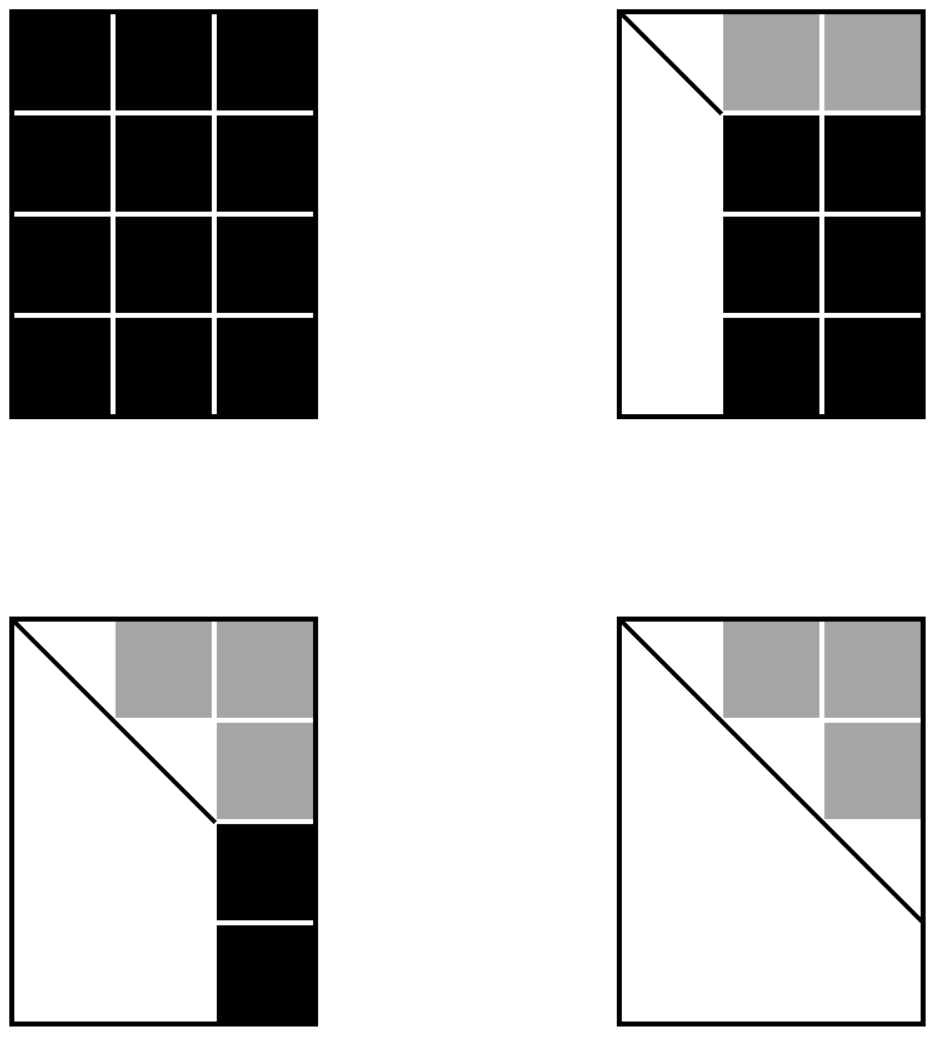}}
\put(006,48){\footnotesize$\mtx{A}^{(0)} = \mtx{A}$}
\put(045,48){\footnotesize$\mtx{A}^{(1)} = (\mtx{U}^{(1)})^{*}\mtx{A}^{(0)}\mtx{V}^{(1)}$}
\put( -2,01){\footnotesize$\mtx{A}^{(2)} = (\mtx{U}^{(2)})^{*}\mtx{A}^{(1)}\mtx{V}^{(2)}$}
\put(045,01){\footnotesize$\mtx{A}^{(3)} = (\mtx{U}^{(3)})^{*}\mtx{A}^{(2)}\mtx{V}^{(3)}$}
\end{picture}
\caption{Cartoon illustrating the process by which a given matrix $\mtx{A}$ is
driven to block upper triangular form,
with most of the mass concentrated into the diagonal entries.
With $\mtx{A}^{(0)} = \mtx{A}$, we form a sequence of matrices
$\mtx{A}^{(i)} = (\mtx{U}^{(i)})^{*}\mtx{A}^{(i-1)}\mtx{V}^{(i)}$ by applying
unitary matrices $\mtx{U}^{(i)}$ and $\mtx{V}^{(i)}$ from the left and the right.
Each $\mtx{U}^{(i)}$ and $\mtx{V}^{(i)}$
consists predominantly of a product of $b$ Householder reflectors, where $b$ is the
block size. Black blocks represent non-zero entries. Grey blocks represent entries that are
not necessarily zero, but are small in magnitude.}
\label{fig:block_cartoon}
\end{figure}

This section describes at a high level the overall template of the algorithm
\randUTV{} that given an $m\times n$ matrix $\mtx{A}$ computes its
UTV decomposition (\ref{eq:defUTVpre}). For simplicity, we assume that $m \geq n$,
that an upper triangular middle factor $\mtx{T}$ is sought, and that we work with
a block size $b$ that evenly divides $n$ so that the matrix $\mtx{A}$ can be
partitioned into $s$ blocks of $b$ columns each; in other words, we assume
that $n = sb$.
The algorithm \randUTV{} iterates over $s$ steps,
where at the $i$'th step the $i$'th block of $b$ columns is driven to upper
triangular form via the application of unitary transformations from the left
and the right.
A cartoon of the process is given in Figure \ref{fig:block_cartoon}.

To be slightly more precise, we build at the $i$'th step unitary matrices $\mtx{U}^{(i)}$
and $\mtx{V}^{(i)}$ of sizes $m\times m$ and $n\times n$, respectively, such that
$$
\mtx{U} = \mtx{U}^{(1)}\mtx{U}^{(2)}\cdots\mtx{U}^{(s)},
\qquad\mbox{and}\qquad
\mtx{V} = \mtx{V}^{(1)}\mtx{V}^{(2)}\cdots\mtx{V}^{(s)}.
$$
Using these matrices, we drive $\mtx{A}$ towards upper triangular form
through a sequence of transformations
\begin{center}
\(
\begin{array}{lcl}
  \mtx{A}^{(0)} &=&\ \mtx{A},\\
  \mtx{A}^{(i)} &=&\ \bigl(\mtx{U}^{(i)}\bigr)^{*}\,\mtx{A}^{(i-1)}\,\mtx{V}^{(i)},\qquad i = 1,\,2,\,3,\,\dots,\,s,\\
  \mtx{T}       &=&\ \mtx{A}^{(s)}.
\end{array}
\)
\end{center}

The objective at step $i$ is to transform the $i$'th diagonal block to diagonal form,
to zero out all blocks beneath the $i$'th diagonal block, and to make all blocks to the
right of the $i$'th diagonal block as small in magnitude as possible.

Each matrix $\mtx{U}^{(i)}$ and $\mtx{V}^{(i)}$ consists predominantly of a product of $b$ Householder
reflectors. To be precise, each such matrix is a product of $b$ Householder reflectors, but with the
$i$'th block of $b$ columns right-multiplied by a small $b\times b$ unitary matrix.

The next two sections provide additional details. Section \ref{sec:singlestep} describes exactly how
to build the transformation matrices $\mtx{U}^{(1)}$ and $\mtx{V}^{(1)}$ that are needed in the
first step of the iteration.
Then, Section \ref{sec:randUTV} shows how to apply the techniques
described in Section \ref{sec:singlestep} repeatedly to build the full factorization.

\section{A randomized algorithm for finding a pair of transformation matrices for the first step}
\label{sec:singlestep}

\subsection{Objectives for the construction}

In this section, we describe a randomized algorithm for finding unitary matrices
$\mtx{U}$ and $\mtx{V}$ that
execute the first step of the process outlined in Section \ref{sec:overview}, and
illustrated in Figure \ref{fig:block_cartoon}. To avoid notational clutter, we
simplify the notation slightly, describing what we consider a basic single step of
the process. Given an $m\times n$ matrix $\mtx{A}$ and a block size $b$, we seek
two orthonormal matrices $\mtx{U}$ and $\mtx{V}$ of sizes $m\times m$ and $n\times n$,
respectively, such that the matrix
\begin{equation}
\label{eq:sia}
\mtx{T} = \mtx{U}^{*}\mtx{A}\mtx{V}
\end{equation}
has a diagonal leading $b\times b$ block, and the entries beneath this block are
all zeroed out. In Section \ref{sec:overview}, we referred to the matrices
$\mtx{U}$ and $\mtx{V}$ as $\mtx{U}^{(1)}$ and $\mtx{V}^{(1)}$, respectively,
and $\mtx{T}$ as $\mtx{A}^{(1)}$.

To make the discussion precise, let us partition $\mtx{U}$ and $\mtx{V}$ so that
\begin{equation}
\label{eq:part1}
\mtx{U} = \bigl[\mtx{U}_{1}\ \ |\ \ \mtx{U}_{2}\bigr],
\qquad\mbox{and}\qquad
\mtx{V} = \bigl[\mtx{V}_{1}\ \ |\ \ \mtx{V}_{2}\bigr],
\end{equation}
where $\mtx{U}_{1}$ and $\mtx{V}_{1}$ each contain $b$ columns.
Then, set
$
\mtx{T}_{ij} = \mtx{U}_{i}^{*}\mtx{A}\mtx{V}_{j}
$
for $i,j = 1,2$ so that
\begin{equation}
\label{eq:part2}
\mtx{U}^{*}\mtx{A}\mtx{V} =
\left[\begin{array}{c|c}\mtx{T}_{11} & \mtx{T}_{12} \\ \hline \mtx{T}_{21} & \mtx{T}_{22}\end{array}\right].
\end{equation}
Our objective is now to build matrices $\mtx{U}$ and $\mtx{V}$ that accomplish the following:
\renewcommand{\theenumi}{\roman{enumi}}
\begin{enumerate}
\item $\mtx{T}_{11}$ is diagonal, with entries that closely approximate the leading $b$ singular values of $\mtx{A}$.
\item $\mtx{T}_{21} = \mtx{0}$.
\item $\mtx{T}_{12}$ has small magnitude.
\item The norm of $\mtx{T}_{22}$ should be close to optimally small, so that $\|\mtx{T}_{22}\| \approx \inf\{\|\mtx{A} - \mtx{C}\|\,\colon\,\mtx{C}\mbox{ has rank }b\}$.
\end{enumerate}
The purpose of condition (iv) is to minimize the error in the rank-$b$ approximation to $\mtx{A}$, cf.~Section \ref{sec:theory}.

\subsection{A theoretically ideal choice of transformation matrices}
\label{sec:optimal}

Suppose that we could somehow find two matrices $\mtx{U}$ and $\mtx{V}$ whose first $b$ columns
\textit{exactly} span the subspaces spanned by the dominant left and right singular vectors, respectively.
Finding such  matrices is of course computationally hard, but \textit{if} we could build them,
then we would  get the optimal result that
\begin{enumerate}
\setcounter{enumi}{1}
\item $\mtx{T}_{21} = \mtx{0}$.
\item $\mtx{T}_{12} = \mtx{0}$.
\item $\|\mtx{T}_{22}\| = \sigma_{b+1} = \min\{\|\mtx{A} - \mtx{C}\|\,\colon\,\mbox{ matrix }\mtx{C}\mbox{ has rank }b\}$.
\end{enumerate}
Enforcing condition (i) is then very easy, since the dominant $b\times b$ block is now
disconnected from the rest of the matrix. Simply executing a full SVD of this small
$b\times b$ block, and then updating $\mtx{U}$ and $\mtx{V}$ will do the job.

\subsection{A randomized technique for approximating the span of the dominant singular vectors}
\label{sec:right}

Inspired by the observation in Section \ref{sec:optimal}
that a theoretically ideal right transformation matrix $\mtx{V}$ is
a matrix whose $b$ first columns span the space spanned by the dominant $b$ right singular vectors of $\mtx{A}$,
we use the randomized power iteration described in Section \ref{sec:randrange} to execute this task.
We let $q$ denote a small integer specifying how many steps of power iteration we take.
Typically, $q=0,1,2$ are good choices.
Then, take the following steps:
(1) Draw a Gaussian random matrix $\mtx{G}$ of size $m \times b$.
(2) Compute a sampling matrix $\mtx{Y} = \bigl(\mtx{A}^{*}\mtx{A}\bigr)^{q}\mtx{A}^{*}\mtx{G}$ of size $n\times b$.
(3) Perform an unpivoted Householder QR factorization of $\mtx{Y}$ so that
$$
\begin{array}{ccccccccccccc}
\mtx{Y} &=& \mathbb{V} & \mtx{R}. \\
n\times b && n\times n & n\times b
\end{array}
$$
Observe that $\mathbb{V}$ will be a product of $b$ Householder reflectors, and that
the first $b$ columns of $\mathbb{V}$ will form an orthonormal basis for the space
spanned by the columns of $\mtx{Y}$. In consequence, the first $b$ columns of
$\mathbb{V}$ form an orthonormal basis for a space that approximately spans
the space spanned by the $b$ dominant right singular vectors of $\mtx{A}$.
(The font used for $\mathbb{V}$ is a reminder
that it is a product of Householder reflectors, which is exploited when it is stored
and operated upon, cf.~Section \ref{sec:tallthin}.)


\subsection{Construction of the left transformation matrix}
\label{sec:left}

The process for finding the left transformation matrix $\mtx{U}$ is deterministic, and will
exactly transform the first $b$ columns of $\mtx{A}\mathbb{V}$ into a diagonal matrix.
Observe that with $\mathbb{V}$ the unitary matrix constructed using the procedure in
Section \ref{sec:right}, we have the identity.
\begin{equation}
\label{eq:def_Aprime}
\mtx{A} =
\mtx{A}\mathbb{V}\mathbb{V}^{*} =
\bigl[\mtx{A}\mathbb{V}_{1}\ \ |\ \ \mtx{A}\mathbb{V}_{2}\bigr]\mathbb{V}^{*},
\end{equation}
where the partitioning $\mathbb{V} = [\mathbb{V}_{1}\ \mathbb{V}_{2}]$ is such that
$\mathbb{V}_{1}$ holds the first $b$ columns of $\mathbb{V}$.
We now perform a full SVD on the matrix $\mtx{A}\mathbb{V}_{1}$, which is of size $m\times b$ so that
\begin{equation}
\label{eq:adele1}
\begin{array}{ccccccccccccc}
\mtx{A}\mathbb{V}_{1} &=& \mtx{U} & \mtx{D} & \mtx{V}_{\rm small}^{*}. \\
m\times b && m\times m & m\times b & b\times b
\end{array}
\end{equation}
Inserting (\ref{eq:adele1}) into (\ref{eq:def_Aprime}) we obtain the identity
\begin{equation}
\label{eq:adele2}
\mtx{A} = \bigl[\mtx{U}\mtx{D}\mtx{V}_{\rm small}^{*}\ \ |\ \ \mtx{A}\mathbb{V}_{2}\bigr]\mathbb{V}^{*}.
\end{equation}
Factor out $\mtx{U}$ in (\ref{eq:adele2}) to the left to get
$$
\mtx{A} = \mtx{U}\bigl[\mtx{D}\mtx{V}_{\rm small}^{*}\ \ |\ \ \mtx{U}^{*}\mtx{A}\mathbb{V}_{2}\bigr]\mathbb{V}^{*}.
$$
Finally, factor out $\mtx{V}_{\rm small}$ to the right to yield the factorization
\begin{equation}
\label{eq:cold}
\mtx{A} = \mtx{U}\underbrace{\bigl[\mtx{D}\ \ |\ \ \mtx{U}^{*}\mtx{A}\mathbb{V}_{2}\bigr]}_{=\mtx{T}}\mtx{V}^{*},
\qquad\mbox{with}\
\mtx{V} = \mathbb{V}\left[\begin{array}{cc}\mtx{V}_{\rm small} & \mtx{0} \\
                                           \mtx{0} & \mtx{I}_{n-b}\end{array}\right].
\end{equation}
Equation (\ref{eq:cold}) is the factorization $\mtx{A} = \mtx{U}\mtx{T}\mtx{V}^{*}$ that we seek, with
$$
\mtx{T} =
\left[\begin{array}{cc}\mtx{T}_{11} & \mtx{T}_{12} \\ \mtx{T}_{21} & \mtx{T}_{22}\end{array}\right] =
\left[\begin{array}{cc}\mtx{D}(1:b,1:b) & \mtx{U}_{1}^{*}\mtx{A}\mathbb{V}_{2} \\ \mtx{0} & \mtx{U}_{2}^{*}\mtx{A}\mathbb{V}_{2}\end{array}\right].
$$

\begin{remark} \label{rem:econ_svd}
As we saw in (\ref{eq:cold}), the right transformation matrix $\mtx{V}$ consists of a product
$\mathbb{V}$ of $b$ Householder reflectors, with the first $b$ columns rotated by a
small unitary matrix $\mtx{V}_{\rm small}$. We will next demonstrate that the left transformation
matrix $\mtx{U}$ can be built in such a way that it can be written in an entirely analogous form.
Simply observe that the matrix
$\mtx{A}\mathbb{V}_{1}$ in (\ref{eq:adele1}) is a tall thin matrix, so that the SVD in (\ref{eq:adele1})
can efficiently be computed by first performing an unpivoted QR factorization
of $\mtx{A}\mathbb{V}_{1}$ to yield a factorization
$$
\begin{array}{ccccccccccccccccc}
\mtx{A}\mathbb{V}_{1} &=& \mathbb{U} & \vtwo{\mtx{R}_{11}}{\mtx{0}},\\
m\times b && m\times m & m\times b
\end{array}
$$
where $\mtx{R}_{11}$ is of size $b\times b$, and $\mathbb{U}$ is a product of $b$
Householder reflectors, cf.~Section \ref{sec:tallthin}.
Then, perform an SVD of $\mtx{R}_{11}$ to obtain
$$
\begin{array}{cccccccccccccccc}
\mtx{R}_{11} &=& \mtx{U}_{\rm small} & \mtx{D}_{\rm small} & \mtx{V}_{\rm small}^{*}.\\
b\times b && b\times b & b\times b & b\times b
\end{array}
$$
The factorization (\ref{eq:adele1}) then becomes
\begin{equation}
\label{eq:warm}
\begin{array}{ccccccccccccc}
\mtx{A}\mathbb{V}_{1} &=& \mtx{U} & \vtwo{\mtx{D}_{\rm small}}{\mtx{0}} & \mtx{V}_{\rm small}^{*}, \\
m\times b && m\times m & m\times b & b\times b
\end{array}
\qquad\mbox{with}\
\mtx{U} = \mathbb{U }\left[\begin{array}{cc} \mtx{U}_{\rm small} & \mtx{0} \\ \mtx{0} & \mtx{I}\end{array}\right].
\end{equation}
We see that the expression for $\mtx{U}$ in (\ref{eq:warm}) is exactly analogous to the expression for
$\mtx{V}$ in (\ref{eq:cold}).
\end{remark}

\subsection{Summary of the construction of transformation matrices}
\label{sec:summary_one_step}

Even though the derivation of the transformation matrices got slightly long, the final algorithm
is simple. It can be written down precisely with just a few lines of Matlab code, as shown
in the right panel in Figure \ref{fig:randUTVmatlab} (the single step process described in
this section is the subroutine \texttt{stepUTV}).

\begin{figure}

\begin{tabular}{|c|c|}\hline
\begin{minipage}{0.52\textwidth}
\small
\begin{verbatim}

function [U,T,V] = randUTV(A,b,q)
  T = A;
  U = eye(size(A,1));
  V = eye(size(A,2));
  for i = 1:ceil(size(A,2)/b)
    I1 = 1:(b*(i-1));
    I2 = (b*(i-1)+1):size(A,1);
    J2 = (b*(i-1)+1):size(A,2);
    if (length(J2) > b)
      [UU,TT,VV] = stepUTV(T(I2,J2),b,q);
    else
      [UU,TT,VV] = svd(T(I2,J2));
    end
    U(:,I2)  = U(:,I2)*UU;
    V(:,J2)  = V(:,J2)*VV;
    T(I2,J2) = TT;
    T(I1,J2) = T(I1,J2)*VV;
  end
return

\end{verbatim}
\end{minipage}
&
\begin{minipage}{0.45\textwidth}
\small
\begin{verbatim}

function [U,T,V] = stepUTV(A,b,q)
  G = randn(size(A,1),b);
  Y = A'*G;
  for i = 1:q
    Y = A'*(A*Y);
  end
  [V,~]    = qr(Y);
  [U,D,W]  = svd(A*V(:,1:b));
  T        = [D,U'*A*V(:,(b+1):end)];
  V(:,1:b) = V(:,1:b)*W;
return









\end{verbatim}
\end{minipage}
\\[2mm] \hline
\end{tabular}
\caption{Matlab code for the algorithm \randUTV{} that given an $m\times n$ matrix $\mtx{A}$
computes its UTV factorization $\mtx{A} = \mtx{U}\mtx{T}\mtx{V}^{*}$, cf.~(\ref{eq:defUTVpre}).
The input parameters $b$ and $q$ reflect the block size and the number of steps of power iteration,
respectively. This code is simplistic in that products of Householder reflectors are stored simply
as dense matrices, making the overall complexity $O(n^{4})$; it also assumes that $m\geq n$. An
efficient implementation is described in Figure \ref{fig:randUTV}.}
\label{fig:randUTVmatlab}
\end{figure}


\section{The algorithm \randUTV{}}
\label{sec:randUTV}

In this section, we describe the algorithm \randUTV{} that given an $m\times n$ matrix $\mtx{A}$
computes a UTV factorization of the form (\ref{eq:defUTVpre}). For concreteness, we assume that
$m\geq n$ and that we seek to build an \textit{upper} triangular middle matrix $\mtx{T}$.
The modifications required for the other cases are straight-forward.
Section \ref{sec:randUTVbasic} describes the most basic version of the scheme,
Section \ref{sec:randUTVeff} describes a computationally efficient version, and
Section \ref{sec:cost} provides a calculation of the asymptotic flop count
of the resulting algorithm.

\subsection{A simplistic algorithm}
\label{sec:randUTVbasic}

The algorithm \randUTV{} is obtained by applying the single-step algorithm
described in Section \ref{sec:singlestep} repeatedly, to drive $\mtx{A}$ to upper
triangular form one block of $b$ columns at a time. We recall that a cartoon of
the process is shown in Figure \ref{fig:block_cartoon}. At the start of the process,
we create three arrays that hold the output matrices $\mtx{T}$, $\mtx{U}$, and $\mtx{V}$,
and initialize them by setting
$$
\mtx{T} = \mtx{A},\qquad
\mtx{U} = \mtx{I}_{m},\qquad
\mtx{V} = \mtx{I}_{n}.
$$
In the first step of the iteration, we use the single-step technique described in
Section \ref{sec:singlestep} to create two unitary ``left and right transformation matrices''
$\mtx{U}^{(1)}$ and $\mtx{V}^{(1)}$ and then update $\mtx{T}$, $\mtx{U}$, and $\mtx{V}$ accordingly:
$$
\mtx{T} \leftarrow (\mtx{U}^{(1)})^{*}\mtx{T}\mtx{V}^{(1)},\qquad
\mtx{U} \leftarrow \mtx{U}\mtx{U}^{(1)},\qquad
\mtx{V} \leftarrow \mtx{V}\mtx{V}^{(1)}.
$$
This leaves us with a matrix $\mtx{T}$ whose $b$ leading columns are upper triangular
(like matrix $\mtx{A}^{(1)}$ in Figure \ref{fig:block_cartoon}).
For the second step, we build transformation matrices
$\mtx{U}^{(2)}$ and $\mtx{V}^{(2)}$ by applying the single-step algorithm described in
Section \ref{sec:singlestep} to the remainder matrix $\mtx{T}((b+1):m,(b+1):n)$, and then
updating $\mtx{T}$, $\mtx{U}$, and $\mtx{V}$ accordingly.

The process then continues to drive one block of $b$ columns at a time to upper triangular
form. With $s = \lceil n/b \rceil$ denoting the total number of steps, we find that after
$s-1$ steps, all that remains to process is the bottom right block of $\mtx{T}$ (cf.~the matrix
$\mtx{A}^{(2)}$ in Figure \ref{fig:block_cartoon}). This block consists of $b$ columns if $n$ is a
multiple of the block size, and is otherwise even thinner. For this last block, we obtain the
final left and right transformation matrices $\mtx{U}^{(s)}$ and $\mtx{V}^{(s)}$ by computing a full
singular value decomposition of the remaining matrix $\mtx{T}(((s-1)b+1):m,((s-1)b+1):n)$, and
updating the matrices $\mtx{T}$, $\mtx{U}$, and $\mtx{V}$ accordingly. (In the cartoon in
Figure \ref{fig:block_cartoon}, we have $s=3$, and the matrices $\mtx{U}^{(3)}$ and $\mtx{V}^{(3)}$
are built by computing a full SVD of a dense matrix of size $2s \times s$.)

We call the algorithm described in this section \randUTV{}.
It can be coded
using just a few lines of Matlab code, as illustrated in Figure \ref{fig:randUTVmatlab}.
In this simplistic version, all unitary matrices are represented as dense matrices, which
makes the overall complexity $O(n^{4})$ for an $n\times n$ matrix.

\subsection{A computationally efficient version}
\label{sec:randUTVeff}

In this section, we describe how the basic version of \randUTV{}, as given in
Figure \ref{fig:randUTVmatlab}, can be turned into a highly efficient procedure via
three simple modifications. The resulting algorithm is summarized in Figure \ref{fig:randUTV}.
We note that the two versions of \randUTV{} described in Figures \ref{fig:randUTVmatlab}
and \ref{fig:randUTV} are mathematically equivalent; if they were to be executed in exact
arithmetic, their outputs would be identical.

The first modification is that all operations on the matrix $\mtx{T}$ and on the
unitary matrices $\mtx{U}$ and $\mtx{V}$ should be carried out ``in place'' to not
unnecessarily move any data. To be precise, using the notation in Section
\ref{sec:singlestep}, we generate at each iteration four unitary matrices
$\mathbb{U}$, $\mathbb{V}$, $\mtx{U}_{\rm small}$, and $\mtx{V}_{\rm small}$. As soon
as such a matrix is generated, it is immediately applied to $\mtx{T}$, and then used to update
either $\mtx{U}$ or $\mtx{V}$.

Second, we exploit that the two ``large'' unitary transforms $\mathbb{U}$ and $\mathbb{V}$ both
consist of products of $b$ Householder reflectors. We generate them by computing unpivoted
Householder QR factorizations of tall and thin matrices, using a subroutine that outputs
simply the $b$ Householder vectors.
Then, $\mathbb{U}$ and $\mathbb{V}$ can both be stored
and applied efficiently, as described in Section \ref{sec:tallthin}.

The third and final modification pertains to the situation where the input matrix $\mtx{A}$ is
non-square. In this case, the full SVD that is computed in the last step involves a
rectangular matrix. When $\mtx{A}$ is tall ($m > n$), we find at this step that $J_{3}$ is
empty, so the matrix to be processed is $\mtx{T}([I_{2},I_{3}],J_{2})$. When computing the
SVD of this matrix, we use the efficient version described in Remark \ref{rem:econ_svd}, which outputs a
factorization in which $\mtx{U}_{\rm small}$ consists in part of a product of Householder
reflectors. (In this situation $\mtx{U}_{\rm small}$ is in fact not necessarily ``small,''
but it can be stored and applied efficiently.)

\begin{remark}[The case $m < n$]
In a situation where the matrix has fewer rows than columns, but we still seek an
upper triangular matrix $\mtx{T}$, \randUTV{} proceeds exactly as described
for the first $s-1$ steps. In the final step, we now find that $I_{3}$ is empty,
but $J_{3}$ is not, and so we need to compute the SVD of the ``fat'' matrix
$\mtx{T}(I_{2},[J_{2},J_{3}])$. We do this in a manner entirely analogous to how
we handle a ``tall'' matrix, by first performing an unpivoted QR factorization
of the \textit{rows} of $\mtx{T}(I_{2},[J_{2},J_{3}])$. In this situation it is
the matrix of right singular vectors at the last step that consists in part of a
product of Householder reflectors.
\end{remark}

\begin{figure}
\fbox{\begin{minipage}{\textwidth}\small
  \begin{algorithmic}

  \STATE $[ \mtx{U}, \, \mtx{T}, \, \mtx{V} ] =
         \randUTV{}( \, \mtx{A}, \, b, \, q )$

  \STATE \textit{\% Initialize output variables:}
  \STATE $\mtx{T} = \mtx{A}$; $\mtx{U} = \mtx{I}_{m}$; $\mtx{V} = \mtx{I}_{n}$;
  \medskip

  \FOR{ $i = 1:\min(\lceil m/b\rceil,\lceil n/b\rceil$) }
    \STATE \textit{\% Create partitions $1:m = [I_{1},I_{2},I_{3}]$ and $1:n = [J_{1},J_{2},J_{3}]$ so that $(I_{2},J_{2})$ points to the ``active'' block} \textit{that is to be diagonalized:}
    \STATE $I_{1} = 1:(b(i-1))$;\ $I_{2} = (b(i-1)+1):\min(bi,m)$;\ $I_{3} = (bi+1):m$;
    \STATE $J_{1} = 1:(b(i-1))$;\ $J_{2} = (b(i-1)+1):\min(bi,n)$;\ $J_{3} = (bi+1):n$;
    \medskip

    \IF{ ($I_{3}$ and $J_{3}$ are both nonempty) }
      \STATE \textit{\% Generate the sampling matrix $\mtx{Y}$ whose columns approximately span the space spanned by the $b$ dominant right singular vectors of the matrix $\mtx{X} = \mtx{T}([I_{2},\,I_{3}],[J_{2},\,J_{3}])$. We do this via randomized sampling, setting $\mtx{Y} = \bigl(\mtx{X}^{*}\mtx{X}\bigr)^{q}\mtx{X}^{*}\mtx{G}$ where $\mtx{G}$ is a Gaussian random matrix with $b$ columns.}
      \STATE $\mtx{G} = \texttt{randn}(m-b(i-1),b)$
      \STATE $\mtx{Y} = \mtx{T}([I_{2},I_{3}],[J_{2},J_{3}])^{*}\mtx{G}$
      \FOR{ $j=1:q$ }
        \STATE $\mtx{Y} = \mtx{T}([I_{2},I_{3}],[J_{2},J_{3}])^{*}\bigl(\mtx{T}([I_{2},I_{3}],[J_{2},J_{3}])\mtx{Y}\bigr)$.
      \ENDFOR
      \medskip

      \STATE \textit{\% Build a unitary matrix $\mathbb{V}$ whose first $b$ columns form an ON basis for the columns of $\mtx{Y}$.  Then, apply the transformations. (We exploit that $\mathbb{V}$ is a product of $b$ Householder reflectors, cf.~Section \ref{sec:tallthin}.)}
      \STATE $[\mathbb{V},\sim] = \texttt{qr}(\mtx{Y})$
      \STATE $\mtx{T}(:,[J_{2},\,J_{3}])  = \mtx{T}(:,[J_{2},J_{3}])\mathbb{V}$
      \STATE $\mtx{V}(:,[J_{2},\,J_{3}])  = \mtx{V}(:,[J_{2},J_{3}])\mathbb{V}$
      \medskip

      \STATE \textit{\% Build a unitary matrix $\mathbb{U}$ whose first $b$ columns form an ON basis for the columns of
                 $\mtx{T}([I_{2},I_{3}],J_{2})$. Then,}
      \STATE \textit{apply the transformations. (We exploit that $\mathbb{U}$ is a product of $b$ Householder reflectors, cf.~Section \ref{sec:tallthin}.)}
      \STATE $[\mathbb{U},\mtx{R}] = \texttt{qr}(\mtx{T}([I_{2},I_{3}],J_{2}))$
      \STATE $\mtx{U}(:,[I_{2},I_{3}])  = \mtx{U}(:,[I_{2},I_{3}])\mathbb{U}$
      \STATE $\mtx{T}([I_{2},I_{3}],J_{3}) = \mathbb{U}^{*}\mtx{T}([I_{2},I_{3}],J_{3})$
      \STATE $\mtx{T}(I_{3},J_{2})      = \mtx{0}$
      \medskip

      \STATE \textit{\% Perform the local SVD that diagonalizes the active diagonal block. Then, apply the transformations.}
      \STATE $[\mtx{U}_{\rm small},\mtx{D}_{\rm small},\mtx{V}_{\rm small}]  = \texttt{svd}(\mtx{R}(1:b,1:b))$
      \STATE $\mtx{T}(I_{2},J_{2})      = \mtx{D}_{\rm small}$
      \STATE $\mtx{T}(I_{2},J_{3})      = \mtx{U}_{\rm small}^{*}\mtx{T}(I_{2},J_{3})$
      \STATE $\mtx{U}( :,I_{2})      = \mtx{U}(:,I_{2})\mtx{U}_{\rm small}$
      \STATE $\mtx{T}(I_{1},J_{2})      = \mtx{T}(I_{1},J_{2})\mtx{V}_{\rm small}$
      \STATE $\mtx{V}( :,J_{2})      = \mtx{V}(:,J_{2})\mtx{V}_{\rm small}$
    \ELSE
      \medskip

      \STATE \textit{\% Perform the local SVD that diagonalizes the last diagonal block. Then, apply the transformations.  If either $I_{3}$ or $J_{3}$ is long, this should be done economically, cf.~Section \ref{sec:randUTVeff}.}
      \STATE $[\mtx{U}_{\rm small},\mtx{D}_{\rm small},\mtx{V}_{\rm small}] = \texttt{svd}(\mtx{T}([I_{2},I_{3}],[J_{2},J_{3}]))$
      \STATE $\mtx{U}(:,[I_{2},I_{3}])       = \mtx{U}(:,[I_{2},I_{3}])\mtx{U}_{\rm small}$
      \STATE $\mtx{V}(:,[J_{2},J_{3}])       = \mtx{V}(:,[J_{2},J_{3}])\mtx{V}_{\rm small}$
      \STATE $\mtx{T}([I_{2},I_{3}],[J_{2},J_{3}]) = \mtx{D}_{\rm small}$
      \STATE $\mtx{T}(I_{1},[J_{2},J_{3}])      = \mtx{T}(I_{1},[J_{2},J_{3}])\mtx{V}_{\rm small}$
    \ENDIF
  \ENDFOR
  \STATE \textbf{return}
  \end{algorithmic}
\end{minipage}}
\caption{The algorithm \randUTV{} that given an $m\times n$ matrix $\mtx{A}$ computes the
UTV factorization $\mtx{A} = \mtx{U}\mtx{T}\mtx{V}^{*}$, cf.~(\ref{eq:defUTVpre}). The input
parameters $b$ and $q$ reflect the block size and the number of steps of power iteration,
respectively. }
\label{fig:randUTV}
\end{figure}

\subsection{Connection between RSVD and \randUTV{}}
\label{sec:theory}

The proposed algorithm \randUTV{} is directly inspired by the Randomized SVD (RSVD) algorithm
described in Remark \ref{remark:RSVD} (as originally described in
\cite{2006_martinsson_random1_orig,2007_martinsson_PNAS,2009_szlam_power} and later elaborated in
\cite{2011_martinsson_random1,2011_martinsson_randomsurvey}). In this section, we explore this
connection in more detail, and demonstrate that the low-rank approximation error that results from
the ``single-step'' UTV-factorization described in Section \ref{sec:singlestep} is identical to
the error produced by the RSVD (with a twist). This means that the detailed error analysis that is available for
the RSVD (see, e.g., \cite{2015_candes_rsvd_bounds,2015_gu_randomized_subspaceiteration,2011_martinsson_randomsurvey})
immediately applies to the procedure described here. To be precise:

\begin{thm}
Let $\mtx{A}$ be an $m\times n$ matrix, let $b$ be an integer denoting step size such that $1 \leq b < \min(m,n)$,
and let $q$ denote a non-negative integer. Let $\mtx{G}$ denote the $m\times b$ Gaussian matrix drawn
in Section \ref{sec:right}, and let $\mtx{U}$, $\mtx{T}$, and $\mtx{V}$ be the factors in the factorization
$\mtx{A} = \mtx{U}\mtx{T}\mtx{V}^{*}$ built in Sections \ref{sec:right} and \ref{sec:left}, partitioned as
in (\ref{eq:part1}) and (\ref{eq:part2}).
\renewcommand{\theenumi}{\alph{enumi}}
\begin{enumerate}
\item Let $\mtx{Y} = \bigl(\mtx{A}^{*}\mtx{A}\bigr)^{q}\mtx{A}^{*}\mtx{G}$ denote a sampling matrix,
and let $\mtx{Q}$ denote an $n\times b$ orthonormal matrix whose columns form a basis for the column
space of $\mtx{Y}$.
Then, the error $\|\mtx{A} - \mtx{A}\mtx{Q}\mtx{Q}^{*}\|$ precisely equals the
error incurred by the RSVD with $q$ steps of power iteration, as analyzed in \cite[Sec.~10]{2011_martinsson_randomsurvey}.
It holds that
\begin{equation}
\label{eq:approx1}
\|\mtx{A} - \mtx{A}\mtx{Q}\mtx{Q}^{*}\| =
\|\mtx{A} - \mtx{U}_{1}\mtx{T}_{11}\mtx{V}_{1}^{*}\| =
\left\|\vtwo{\mtx{T}_{12}}{\mtx{T}_{22}}\right\|.
\end{equation}
\item Let $\mtx{Z} = \mtx{A}\mtx{Y} =
\bigl(\mtx{A}\mtx{A}^{*}\bigr)^{q+1}\mtx{G}$ denote a sampling matrix,
and let $\mtx{W}$ denote an $m\times b$ orthonormal matrix whose columns form a basis for the column
space of $\mtx{Z}$. If the rank of $\mtx{A}$ is at least $b$, then
\begin{equation}
\label{eq:approx2}
\|\mtx{A} - \mtx{W}\mtx{W}^{*}\mtx{A}\| =
\|\mtx{A} - \mtx{U}_{1}\bigl(\mtx{T}_{11}\mtx{V}_{1}^{*} + \mtx{T}_{12}\mtx{V}_{2}^{*}\bigr)\| =
\left\|\mtx{T}_{22}\right\|.
\end{equation}
\end{enumerate}
\end{thm}

We observe that the term $\|\mtx{A} - \mtx{W}\mtx{W}^{*}\mtx{A}\|$ that arises in part (b) can
informally be said to be the error resulting from RSVD with ``$q+1/2$'' steps of power iteration.
This conforms with what one might have optimistically hoped for, given that the RSVD involves
$2q+2$ applications of either $\mtx{A}$ or $\mtx{A}^{*}$ to thin matrices with $b$ columns,
and \randUTV{} involves $2q+3$ such operations at each step ($2q+1$ applications in
building $\mtx{Y}$, and then the computations of $\mtx{A}\mathbb{V}$ and $\mathbb{U}^{*}\mtx{A}$,
which are in practice applications of $\mtx{A}$ to thin matrices due to the identity (\ref{eq:WYrepresentation})).

\begin{proof} The proofs for the two parts rest on the fact that $\mtx{A}$ can be decomposed as follows:
\begin{equation}
\label{eq:scoot}
\mtx{A} =
\mtx{U}\mtx{T}\mtx{V}^{*} =
\bigl[\mtx{U}_{1}\ \mtx{U}_{2}\bigr]
\left[\begin{array}{cc}
\mtx{T}_{11} & \mtx{T}_{12} \\
\mtx{0}      & \mtx{T}_{22}
\end{array}\right]
\left[\begin{array}{c}
\mtx{V}_{1}^{*} \\
\mtx{V}_{2}^{*}
\end{array}\right] =
\mtx{U}_{1}\mtx{T}_{11}\mtx{V}_{1}^{*} +
\mtx{U}_{1}\mtx{T}_{12}\mtx{V}_{2}^{*} +
\mtx{U}_{2}\mtx{T}_{22}\mtx{V}_{2}^{*},
\end{equation}
where $\mtx{U}_{1}$ and $\mtx{V}_{1}$ have $b$ columns each, and $\mtx{T}_{11}$ is of size $b\times b$.
With the identity (\ref{eq:scoot}) in hand, the claims in (a) follow once we have established that
the orthogonal projection $\mtx{Q}\mtx{Q}^{*}$ equals the projection $\mtx{V}_{1}\mtx{V}_{1}^{*}$. The
claims in (b) follow once we establish that $\mtx{W}\mtx{W}^{*} = \mtx{U}_{1}\mtx{U}_{1}^{*}$.

\lsp

(a) Observe that the matrix $\mtx{Q}$ is by construction identical to the matrix $\mathbb{V}_{1}$
built in Sections \ref{sec:right} and \ref{sec:left}. Since $\mathbb{V}_{1} = \mtx{V}_{1}\mtx{V}_{\rm small}^{*}$,
we find that
$\mtx{Q}\mtx{Q}^{*} =
\mathbb{V}_{1}\mathbb{V}_{1}^{*} =
\bigl(\mtx{V}_{1}\mtx{V}_{\rm small}^{*}\bigr)
\bigl(\mtx{V}_{1}\mtx{V}_{\rm small}^{*}\bigr)^{*} =
\mtx{V}_{1}\bigl(\mtx{V}_{\rm small}^{*}\mtx{V}_{\rm small}\bigr)\mtx{V}_{1}^{*}$.
Since $\mtx{V}_{\rm small}^{*}\mtx{V}_{\rm small} = \mtx{I}_{b}$, it follows that
$\mtx{Q}\mtx{Q}^{*} = \mtx{V}_{1}\mtx{V}_{1}^{*}$.
Then
\begin{equation}
\label{eq:WD4}
\mtx{A}\mtx{Q}\mtx{Q}^{*} =
\mtx{A}\mtx{V}_{1}\mtx{V}_{1}^{*} = \{\mbox{Use}\ (\ref{eq:scoot})\ \mbox{and that}\ \mtx{V}_{2}^{*}\mtx{V}_{1} = \mtx{0}\ \mbox{and}\ \mtx{V}_{1}^{*}\mtx{V}_{1} = \mtx{I}.\} =
\mtx{U}_{1}\mtx{T}_{11}\mtx{V}_{1}^{*}.
\end{equation}
The first identity in (\ref{eq:approx1}) follows immediately from (\ref{eq:WD4}).
The second identity holds since (\ref{eq:scoot}) implies that
$\mtx{A} - \mtx{U}_{1}\mtx{T}_{11}\mtx{V}_{1}^{*} =
\mtx{U}_{1}\mtx{T}_{12}\mtx{V}_{2}^{*} +
\mtx{U}_{2}\mtx{T}_{22}\mtx{V}_{2}^{*} =
\mtx{U}\vtwo{\mtx{T}_{12}}{\mtx{T}_{22}}\mtx{V}_{2}^{*}$,
with $\mtx{U}$ unitary and $\mtx{V}_{2}$ orthonormal.

\lsp

(b) We will first prove that with probability 1, the two $m\times b$ matrices $\mtx{A}\mathbb{V}_{1}$ and $\mtx{Z}$
have the same column spaces. To this end, note that since $\mathbb{V}_{1}$ is obtained by performing an unpivoted
QR factorization of the matrix $\mtx{Y}$ defined in (a), we know that $\mathbb{V}_{1}\mtx{R} = \mtx{Y}$ for some
$b\times b$ upper triangular matrix $\mtx{R}$. The assumption that $\mtx{A}$ has rank at least $b$ implies that
$\mtx{R}$ is invertible with probability 1. Consequently, $\mtx{A}\mathbb{V}_{1} = \mtx{A}\mtx{Y}\mtx{R}^{-1} =
\mtx{Z}\mtx{R}^{-1}$, since $\mtx{Z} = \mtx{A}\mtx{Y}$. Since right multiplication by an invertible matrix does
not change the column space of a matrix, the claim follows.

\lsp

Since $\mtx{A}\mathbb{V}_{1}$ and $\mtx{Z}$ have the same column spaces, it follows from the definition of $\mathbb{U}$ that
$\mathbb{U}_{1}\mathbb{U}_{1}^{*} = \mtx{W}\mtx{W}^{*}$. Since $\mtx{U}_{1} = \mathbb{U}_{1}\mtx{U}_{\rm small}$ where
$\mtx{U}_{\rm small}$ is unitary, we see that $\mathbb{U}_{1}\mathbb{U}_{1}^{*} = \mtx{U}_{1}\mtx{U}_{1}^{*}$. Consequently,
$$
\mtx{W}\mtx{W}^{*}\mtx{A} =
\mtx{U}_{1}\mtx{U}_{1}^{*}\mtx{A} = \{\mbox{Use}\ (\ref{eq:scoot}).\} =
\mtx{U}_{1}\mtx{T}_{11}\mtx{V}_{1}^{*} +
\mtx{U}_{1}\mtx{T}_{12}\mtx{V}_{2}^{*}.
$$
which establishes the first identity in (\ref{eq:approx2}).
The second identity holds since (\ref{eq:scoot}) implies that
$\mtx{A} - \mtx{U}_{1}(\mtx{T}_{11}\mtx{V}_{1}^{*} + \mtx{T}_{12}\mtx{V}_{2}^{*})
= \mtx{U}_{2}\mtx{T}_{22}\mtx{V}_{2}^{*}$, with $\mtx{U}_{2}$ and $\mtx{V}_{2}$ orthonormal.
\end{proof}

\begin{remark}[Oversampling]
\label{remark:positivep}
We recall that the accuracy of \randUTV{} depends on how well the space
$\mbox{col}(\mathbb{V}_{1})$ aligns with the space spanned by the $b$ dominant
right singular vectors of $\mtx{A}$. If these two spaces were to match exactly,
then the truncated UTV factorization would achieve perfectly optimal accuracy.
One way to improve the alignment is to increase the power parameter $q$. A second
way to make the two spaces align better is to use \textit{oversampling}, as described
in Remark \ref{remark:oversampling}. With $p$ an over-sampling parameter (say $p=5$
or $p=10$), we would draw a Gaussian random matrix $\mtx{G}$ of size $m\times (b+p)$,
and then compute an ``extended'' sampling matrix $\mtx{Y}' = \bigl(\mtx{A}^{*}\mtx{A}\bigr)^{q}\mtx{A}^{*}\mtx{G}$
of size $n\times (b+p)$. The $n\times b$ sampling matrix $\mtx{Y}$ we would use to
compute $\mathbb{V}$ would then be formed by the $b$ dominant left singular vectors
of $\mtx{Y}'$, cf.~Figure \ref{fig:stepUTVoversampling}. Oversampling in this fashion does improve the accuracy (see
Section \ref{app:oversampling}), but in our experience, the
additional computational cost is not worth it. Incorporating over-sampling
would also introduce an additional tuning parameter $p$, which is in many ways
undesirable.
\end{remark}

\subsection{Theoretical cost of \randUTV{}}
\label{sec:cost}

Now we analyze the theoretical cost of the implementation of \randUTV{},
and we compare it to those of CPQR and SVD.

The theoretical cost of the CPQR factorization and
the unpivoted QR factorization of an $m \times n$ matrix is:
$2mn^2 - 2n^3/3$ flops,
when no orthonormal matrices are required.
Although the theoretical cost of both the pivoted and the unpivoted QR
is the same, other factors should be considered,
being the most important one the quality of flops.
In modern architectures,
flops performed inside BLAS-3 operations can be about 5--10 times faster than
flops performed inside BLAS-1 and BLAS-2 operations,
since BLAS-3 is CPU-bound whereas BLAS-1 and BLAS-2 are memory-bound.
Hence, the unpivoted QR factorization in
high-performance libraries such as~\cite{LAPACK3} can be much faster
because most of the flops are performed inside BLAS-3 operations,
whereas only half of the flops
in the best implementations of CPQR (\texttt{dgeqp3}) in \cite{LAPACK3}
are performed inside BLAS-3 operations.
In addition, this low performance of CPQR can even be smaller
because of the appearance of catastrophic cancellations
during the computations.
The appearance of just one catastrophic cancellation
will stop the building of a block Householder reflector
before all of it has been built.
This sudden stop forces the algorithm to work on smaller block sizes,
which are suboptimal, and hence performances are even lower.

The SVD usually comprises two steps:
the reduction to bidiagonal form,
and then the reduction from bidiagonal to diagonal form.
The first step is a direct step, whereas the second step is an iterative one.
The theoretical cost of the reduction to bidiagonal form of
an $m \times n$ matrix is:
$4mn^2 - 4n^3/3$ flops,
when no singular vectors are needed.
If $m \gg n$, the cost can be reduced to:
$2mn^2 + 2n^3$ flops by performing first a QR factorization.
The cost of the reduction to diagonal form depends
on the number of iterations, which is unknown a priori,
but it is usually small when no singular vectors are built.
On the one hand,
in the bidiagonalization a large share of the flops are performed
inside the not-so-efficient BLAS-1 and BLAS-2 operations.
Therefore, no high performances are obtained
in the reduction to bidiagonal form.
On the other hand,
the reduction to diagonal form uses just BLAS-1 operations.
These operations are memory-bound, and
in addition they cannot be efficiently parallelized within BLAS,
which might reduce performances on multicore machines.
In conclusion, usual implementations of the SVD will render
low performances on both single-core architectures and
multicore architectures.

The theoretical cost of the \randUTV{} factorization
of an $m \times n$ matrix is:
$(5 + 2q) mn^2 - (3 + 2q)n^3/3$ flops,
when no orthonormal matrices are required,
and $q$ steps of power iteration are applied.
If $q=0$, the theoretical cost of \randUTV{} is three times as high as
the theoretical cost of CPQR;
if $q=1$, it is four times as high;
and
if $q=2$, it is five times as high.
Although \randUTV{} seems computationally more expensive than CPQR,
the quality of flops should be considered.
The share of BLAS-1 and BLAS-2 flops in \randUTV{} is very small:
BLAS-1 and BLAS-2 flops are only employed inside
the CPQR factorization of the sampling matrix $\mtx{Y}$,
the QR factorization of the current column block,
and the SVD of the diagonal block.
As these operations only apply to blocks of dimensions
$n \times b$, $m \times b$, and $b \times b$, respectively,
the total amount of these types of flops is negligible,
and therefore most of the flops performed in \randUTV{} are BLAS-3 flops.
Hence, the algorithm for computing the \randUTV{} will be much faster than
what the theoretical cost predicts.
In conclusion,
this heavy use of BLAS-3 operations will render
good performances on single-core architectures,
multicore architectures, GPUs,
and distributed-memory architectures.

\section{Numerical results}
\label{sec:num}

\subsection{Computational speed}
\label{sec:speed}

In this section, we investigate the speed of the proposed algorithm \randUTV{},
and compare it to the speeds of highly optimized methods
for computing the SVD and the column pivoted QR (CPQR) factorization.

All experiments reported in this article were performed on an
Intel Xeon E5-2695 v3 (Haswell) processor (2.3 GHz), with 14 cores.
In order to be able to show scalability results, the clock speed was
throttled at 2.3 GHz, turning off so-called turbo boost.
Other details of interest include
that the OS used was Linux (Version 2.6.32-504.el6.x86\_64),
and the code was compiled with gcc (Version 4.4.7).
Main routines for computing the SVD ({\tt dgesvd})
and the CPQR ({\tt dgeqp3}) were taken from
Intel's MKL library (Version 11.2.3)
since this library usually delivers much higher performances
than Netlib's LAPACK codes.
Our implementations were coded
with {\tt libflame}~\cite{CiSE09,libflame_ref} (Release 11104).

Each of the three algorithms we tested (randUTV, SVD, CPQR) was applied
to double-precision real matrices of
size $n\times n$. We report the following times:
\begin{tabbing}
\hspace{5mm}\=\hspace{15mm} \=\hspace{60mm} \= \kill
\> $T_{\rm svd}$     \>
The time in seconds for the LAPACK function \texttt{dgesvd} from Intel's MKL.\\
\> $T_{\rm cpqr}$    \>
The time in seconds for the LAPACK function \texttt{dgeqp3} from Intel's MKL.\\
\> $T_{\rm randUTV}$ \>
The time in seconds for our implementation of \randUTV{}.
\end{tabbing}
For the purpose of a fair comparison,
the three implementations were linked to the BLAS library from Intel's MKL.
In all cases, we used an algorithmic block size of $b = 64$.
While likely not optimal for all problem sizes,
this block size yields near best performance and,
regardless, it allows us to easily compare and
contrast the performance of the different implementations.

\begin{figure}
\begin{center}
\begin{tabular}{cc}
\includegraphics[width=0.40\textwidth]{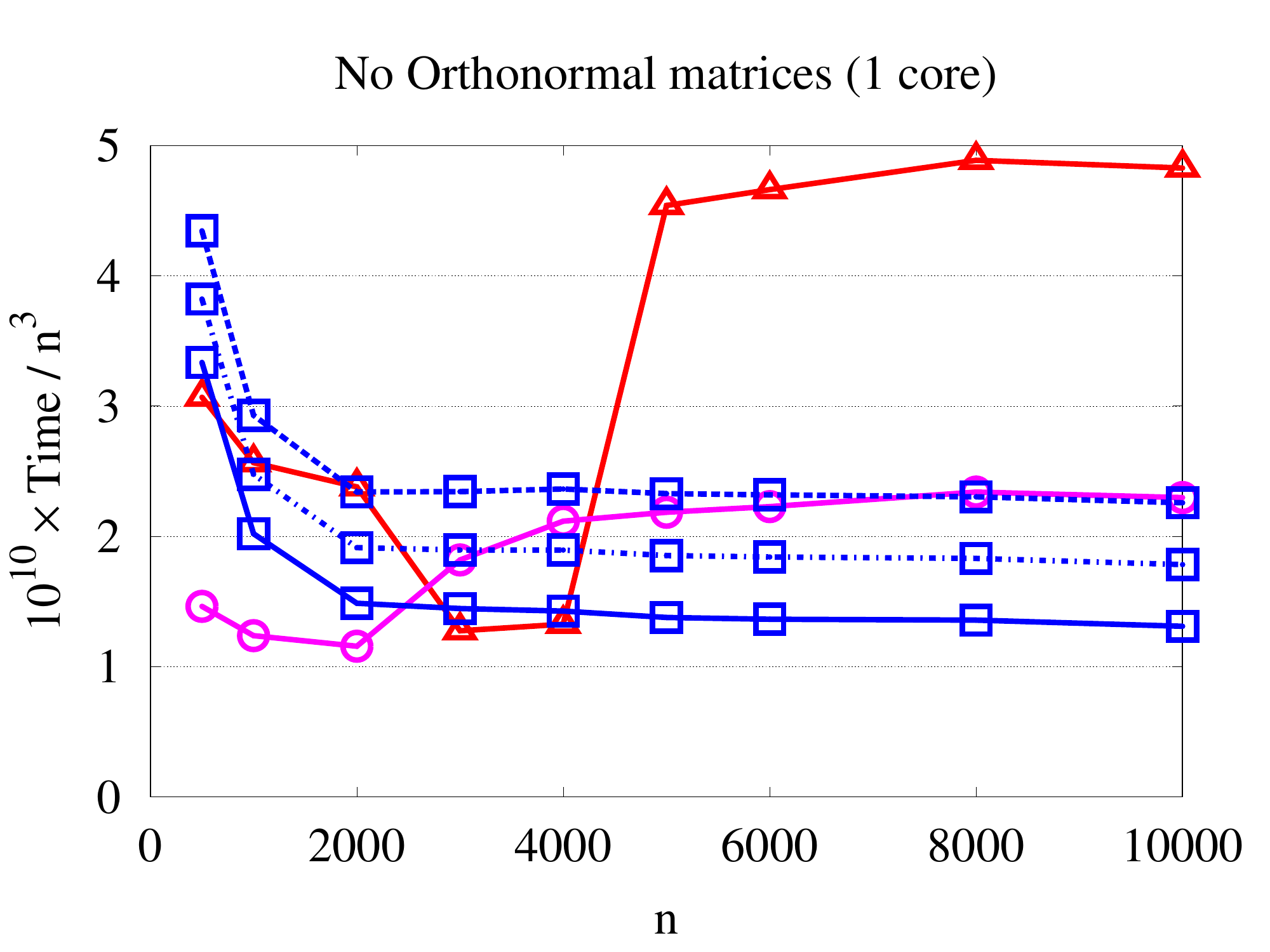} &
\includegraphics[width=0.40\textwidth]{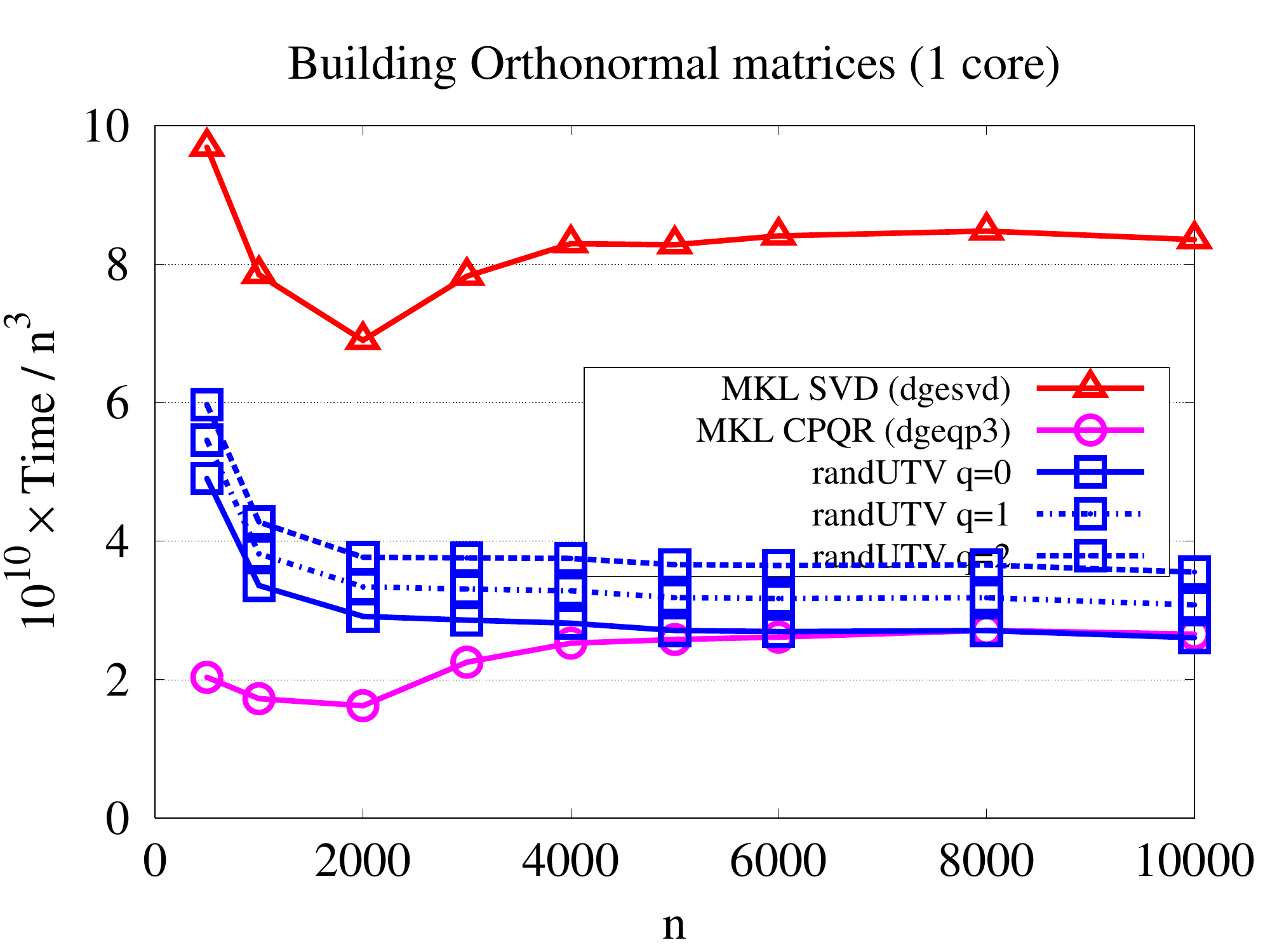} \\
\includegraphics[width=0.40\textwidth]{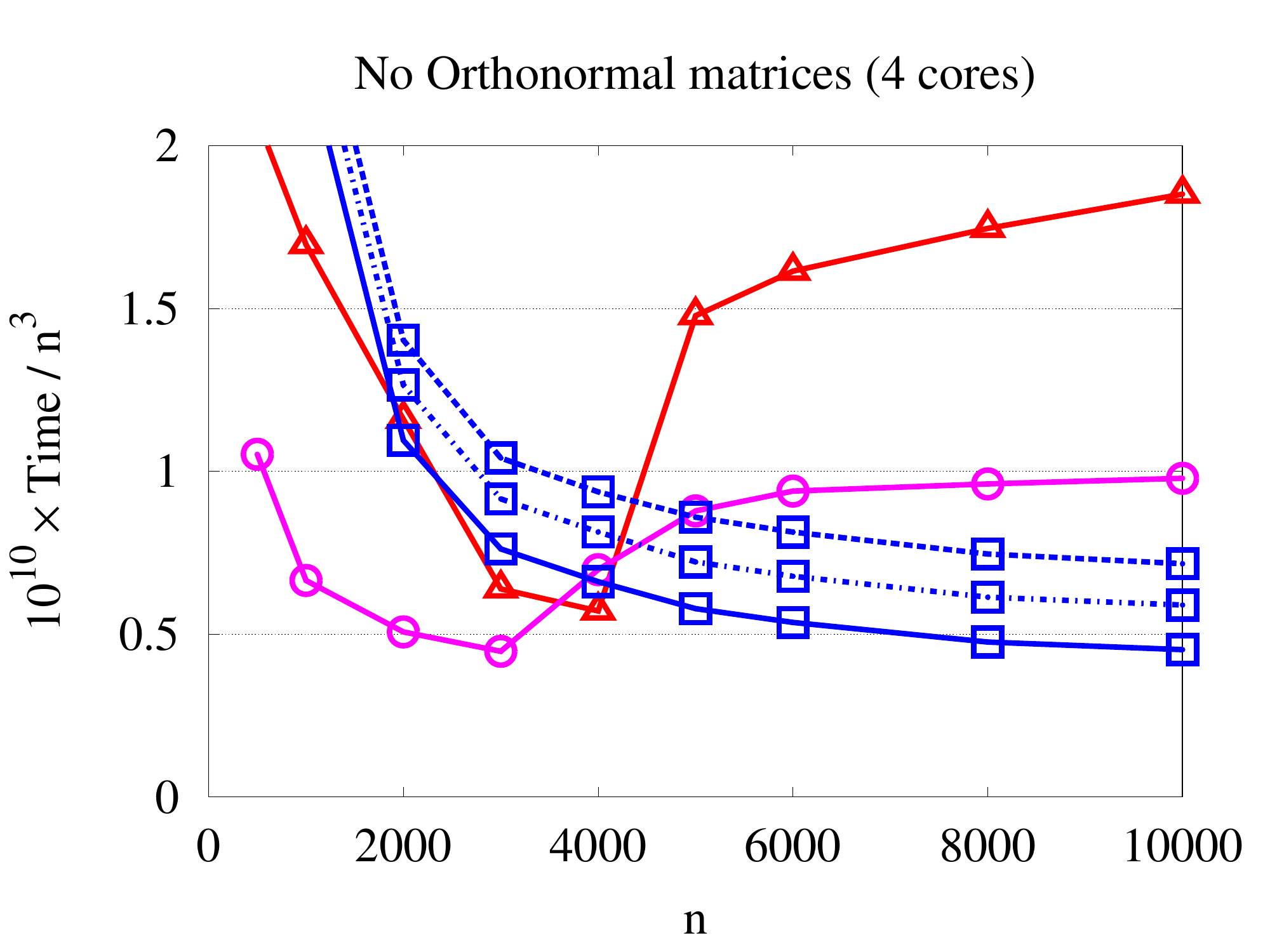} &
\includegraphics[width=0.40\textwidth]{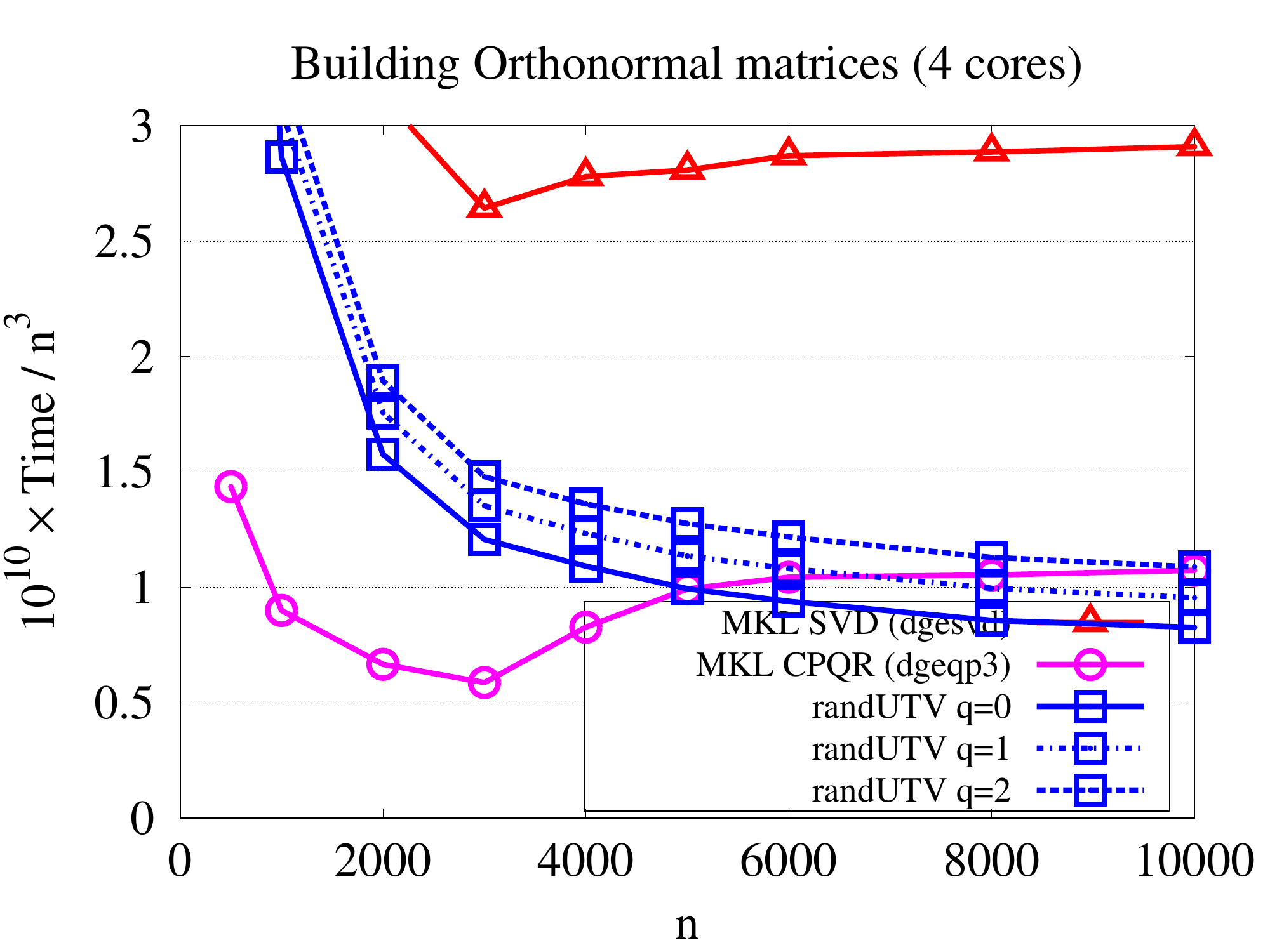} \\
\includegraphics[width=0.40\textwidth]{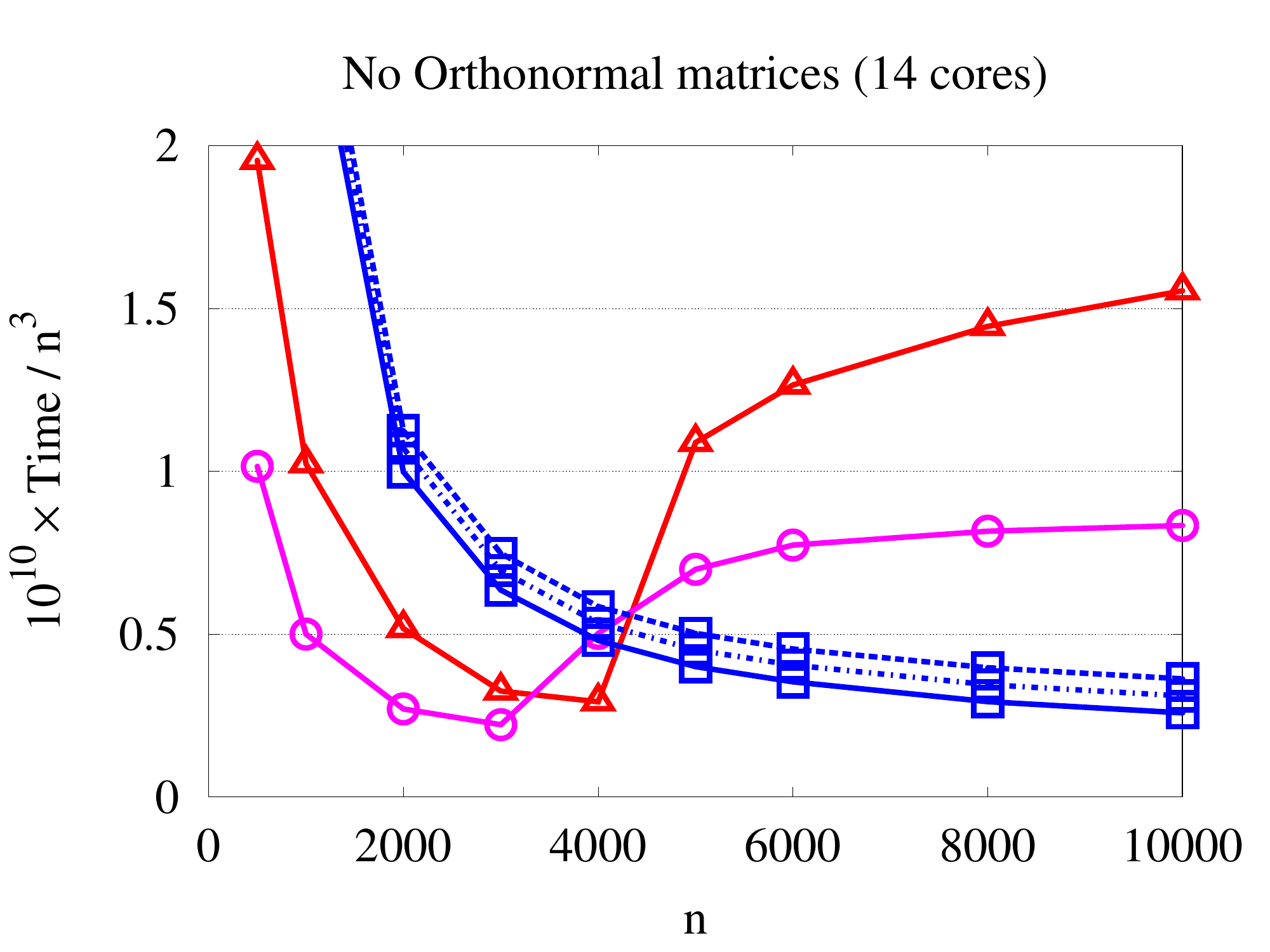} &
\includegraphics[width=0.40\textwidth]{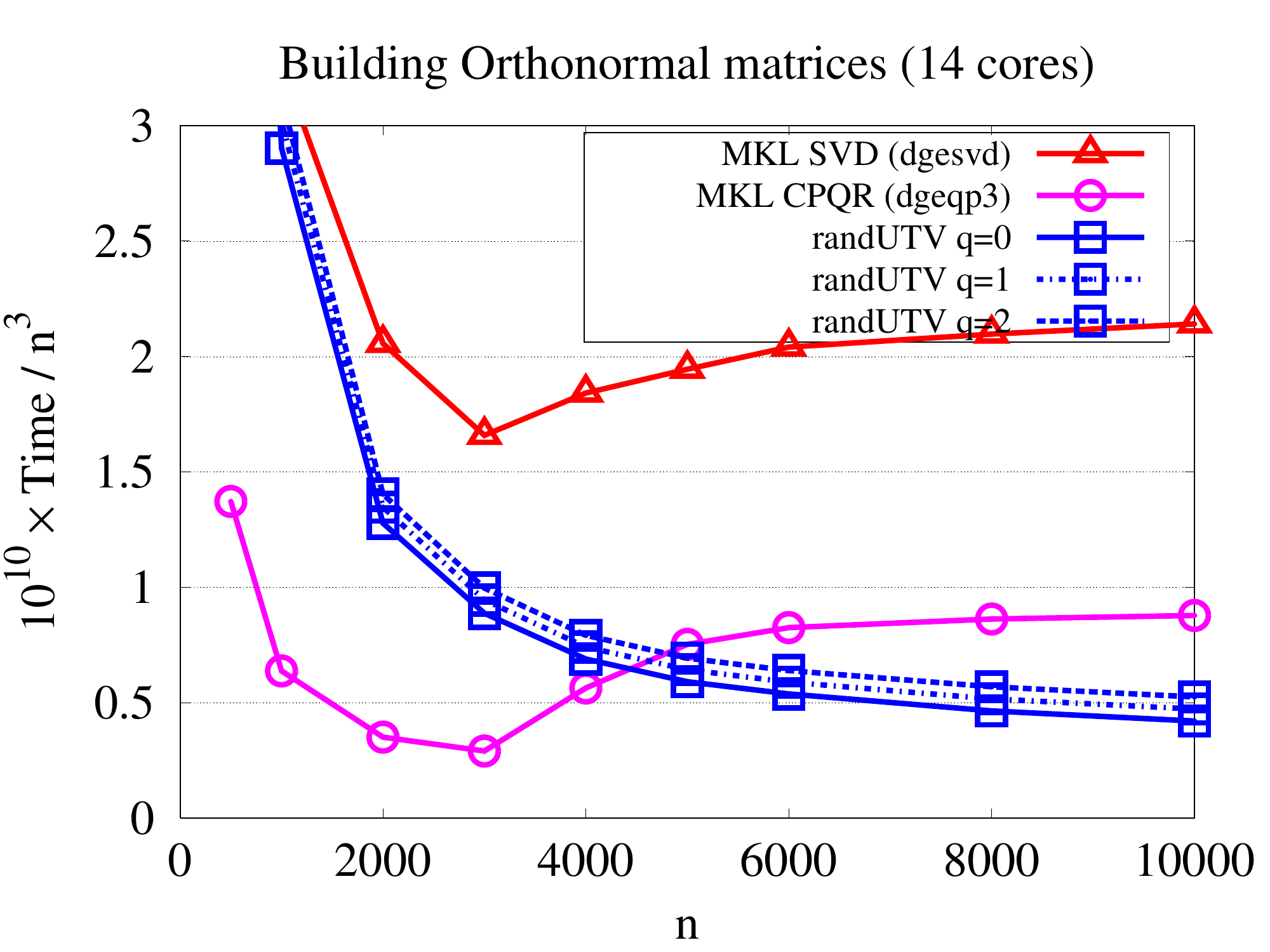} \\
\end{tabular}
\end{center}

\caption{Computational cost (the lower, the better performances) of \randUTV{}
compared to the cost of the LAPACK drivers \texttt{dgesvd} (SVD) and
\texttt{dgeqp3} (CPQR).
The algorithms were applied to double-precision real matrices
of size $n\times n$.
The measured computational time divided by $n^{3}$ is plotted against $n$
(e.g.~the red line is $T_{\rm svd}/n^{3}$.)
The three rows of plots correspond to executing on 1, 4, and 14 cores,
respectively.}
\label{fig:cputimes}
\end{figure}

Table~\ref{table:times_ort} 
shows the measured computational times when executed on 1, 4, and 14 cores,
respectively.
In these experiments,
all orthonormal matrices ($\mtx{U}$ and $\mtx{V}$ for SVD and UTV,
and $\mtx{Q}$ for CPQR) are explicitly formed.
This slightly favors CPQR since only one orthonormal matrix is required.
The corresponding numbers obtained
when orthonormal matrices are not built are given
in Appendix \ref{app:moretimes}.
To better illustrate the relative performance of the various
techniques, we plot in Figure \ref{fig:cputimes}
the computational times measured divided by $n^{3}$.
Since all techniques under consideration have
asymptotic complexity $O(n^{3})$ when applied to an $n\times n$ matrix, these graphs better
reveal the computational efficiency. (We plot time divided by $n^{3}$ rather than the more
commonly reported ``normalized Gigaflops'' since the algorithms we compare have different
scaling factors multiplying the dominant $n^{3}$-term in the asymptotic flop count.) Figure
\ref{fig:cputimes} also shows the timings we measured when the orthonormal matrices were not formed.

The results in Figure \ref{fig:cputimes} lead us to make several observations: 
(1) The algorithm \randUTV{} is decisively faster than the SVD in almost all
cases (the exceptions involve the situation where no unitary matrices are
sought, and the input matrix is small).
(2) Comparing the speeds of CPQR and \randUTV{}, we see that when both methods
are executed on a single core, the speeds are similar, with CPQR being slightly 
faster in some regimes. 
(3) As the matrix size grows, and as the number of cores increases, \randUTV{} 
gains an edge on CPQR in terms of speed.

\begin{table}
\begin{tabular}{r|r|r|rrr}
\multicolumn{6}{c}{Computational times when executed on a single core} \\
\multicolumn{1}{c|}{$n$} &
\multicolumn{1}{c|}{$T_{\rm svd}$} &
\multicolumn{1}{c|}{$T_{\rm cpqr}$} &
\multicolumn{3}{c}{$T_{\rm randUTV}$} \\
\multicolumn{1}{c|}{} &
\multicolumn{1}{c|}{} &
\multicolumn{1}{c|}{} &
\multicolumn{1}{|c}{$q=0$} &
\multicolumn{1}{c}{$q=1$} &
\multicolumn{1}{c}{$q=2$} \\ \hline
    500 &  1.21e-01 &  2.54e-02 &  6.13e-02 &  6.82e-02 &  7.46e-02 \\
   1000 &  7.85e-01 &  1.72e-01 &  3.36e-01 &  3.82e-01 &  4.27e-01 \\
   2000 &  5.52e+00 &  1.30e+00 &  2.33e+00 &  2.67e+00 &  3.01e+00 \\
   3000 &  2.11e+01 &  6.08e+00 &  7.72e+00 &  8.93e+00 &  1.01e+01 \\
   4000 &  5.31e+01 &  1.62e+01 &  1.80e+01 &  2.10e+01 &  2.40e+01 \\
   5000 &  1.04e+02 &  3.22e+01 &  3.39e+01 &  3.98e+01 &  4.57e+01 \\
   6000 &  1.82e+02 &  5.65e+01 &  5.82e+01 &  6.85e+01 &  7.88e+01 \\
   8000 &  4.34e+02 &  1.39e+02 &  1.39e+02 &  1.63e+02 &  1.87e+02 \\
  10000 &  8.36e+02 &  2.66e+02 &  2.60e+02 &  3.08e+02 &  3.55e+02 \\
\end{tabular}
\vspace*{0.5cm}

\begin{tabular}{r|r|r|rrr}
\multicolumn{6}{c}{Computational times when executed on 4 cores} \\
\multicolumn{1}{c|}{$n$} &
\multicolumn{1}{c|}{$T_{\rm svd}$} &
\multicolumn{1}{c|}{$T_{\rm cpqr}$} &
\multicolumn{3}{c}{$T_{\rm randUTV}$} \\
\multicolumn{1}{c|}{} &
\multicolumn{1}{c|}{} &
\multicolumn{1}{c|}{} &
\multicolumn{1}{|c}{$q=0$} &
\multicolumn{1}{c}{$q=1$} &
\multicolumn{1}{c}{$q=2$} \\ \hline
    500 &  8.73e-02 &  1.80e-02 &  7.59e-02 &  7.98e-02 &  8.16e-02 \\
   1000 &  4.20e-01 &  9.00e-02 &  2.86e-01 &  3.07e-01 &  3.22e-01 \\
   2000 &  2.50e+00 &  5.33e-01 &  1.26e+00 &  1.40e+00 &  1.52e+00 \\
   3000 &  7.14e+00 &  1.58e+00 &  3.26e+00 &  3.65e+00 &  3.99e+00 \\
   4000 &  1.78e+01 &  5.29e+00 &  6.99e+00 &  7.89e+00 &  8.71e+00 \\
   5000 &  3.51e+01 &  1.24e+01 &  1.24e+01 &  1.42e+01 &  1.59e+01 \\
   6000 &  6.20e+01 &  2.25e+01 &  2.03e+01 &  2.33e+01 &  2.63e+01 \\
   8000 &  1.48e+02 &  5.39e+01 &  4.39e+01 &  5.10e+01 &  5.78e+01 \\
  10000 &  2.91e+02 &  1.07e+02 &  8.26e+01 &  9.55e+01 &  1.09e+02 \\
\end{tabular}
\vspace*{0.5cm}

\begin{tabular}{r|r|r|rrr}
\multicolumn{6}{c}{Computational times when executed on 14 cores} \\
\multicolumn{1}{c|}{$n$} &
\multicolumn{1}{c|}{$T_{\rm svd}$} &
\multicolumn{1}{c|}{$T_{\rm cpqr}$} &
\multicolumn{3}{c}{$T_{\rm randUTV}$} \\
\multicolumn{1}{c|}{} &
\multicolumn{1}{c|}{} &
\multicolumn{1}{c|}{} &
\multicolumn{1}{|c}{$q=0$} &
\multicolumn{1}{c}{$q=1$} &
\multicolumn{1}{c}{$q=2$} \\ \hline
    500 &  7.51e-02 &  1.71e-02 &  8.90e-02 &  9.18e-02 &  9.35e-02 \\
   1000 &  3.26e-01 &  6.37e-02 &  2.90e-01 &  3.02e-01 &  3.13e-01 \\
   2000 &  1.65e+00 &  2.80e-01 &  1.02e+00 &  1.08e+00 &  1.12e+00 \\
   3000 &  4.48e+00 &  7.83e-01 &  2.40e+00 &  2.55e+00 &  2.69e+00 \\
   4000 &  1.18e+01 &  3.61e+00 &  4.41e+00 &  4.74e+00 &  5.07e+00 \\
   5000 &  2.43e+01 &  9.41e+00 &  7.38e+00 &  8.03e+00 &  8.66e+00 \\
   6000 &  4.41e+01 &  1.78e+01 &  1.16e+01 &  1.27e+01 &  1.38e+01 \\
   8000 &  1.07e+02 &  4.41e+01 &  2.38e+01 &  2.64e+01 &  2.91e+01 \\
  10000 &  2.14e+02 &  8.77e+01 &  4.20e+01 &  4.72e+01 &  5.25e+01 \\
\end{tabular}
\vspace*{0.2cm}
\caption{Computational times of different factorizations
when executed on one core (top), 4 cores (middle), and 14 cores (bottom).
All orthonormal matrices are built explicitly.}
\label{table:times_ort}
\end{table}

\subsection{Errors}
\label{sec:errors}

In this section, we describe the results of the numerical experiments
that were conducted to investigate how good \randUTV{} is
at computing an accurate rank-$k$ approximation to a given matrix.
Specifically,
we compared how well partial factorizations reveal the numerical ranks of four
different test matrices:
\begin{itemize}
\item \textit{Matrix 1 (Fast Decay):} This is an $n\times n$ matrix of the form
$\mtx{A} = \mtx{U}\mtx{D}\mtx{V}^{*}$ where $\mtx{U}$ and $\mtx{V}$ are randomly
drawn matrices with orthonormal columns (obtained by performing an unpivoted QR
factorization on a random Gaussian matrix), and where $\mtx{D}$ is diagonal with entries
given by $d_{j} = \beta^{(j-1)/(n-1)}$ with $\beta = 10^{-5}$.
\item \textit{Matrix 2 (S-shaped Decay):} This matrix is built in the
same manner as ``Matrix 1'', but now the diagonal entries
of $\mtx{D}$ are chosen to first hover around 1, then decay rapidly, and then
level out at $10^{-2}$, as shown in Figure \ref{fig:errors_S_spec} (black line).
\item \textit{Matrix 3 (Gap):} This matrix is built in the
same manner as ``Matrix 1'', but now there is a sudden drop in magnitudes of the singular
values so that $\sigma_{151} = 0.1\,\sigma_{150}$. Specifically, $\mtx{D}(j,j) = 1/j$ for $j \leq 150$,
and $\mtx{D}(j,j) = 0.1/j$ for $j > 150$.
\item \textit{Matrix 4 (BIE):} This matrix is the result of discretizing
a Boundary Integral Equation (BIE) defined on a smooth closed curve in the plane. To be precise,
we discretized the so called ``single layer'' operator associated with the Laplace equation
using a $6^{\rm th}$ order quadrature rule designed by Alpert \cite{1999_alpert_hybrid}.
This operator is well-known to be ill-conditioned, which necessitates the use of a rank-revealing
factorization in order to solve the corresponding linear system in as stable a manner as possible.
\end{itemize}

\begin{figure}
\vspace*{0.3cm}
\includegraphics[width=0.95\textwidth]{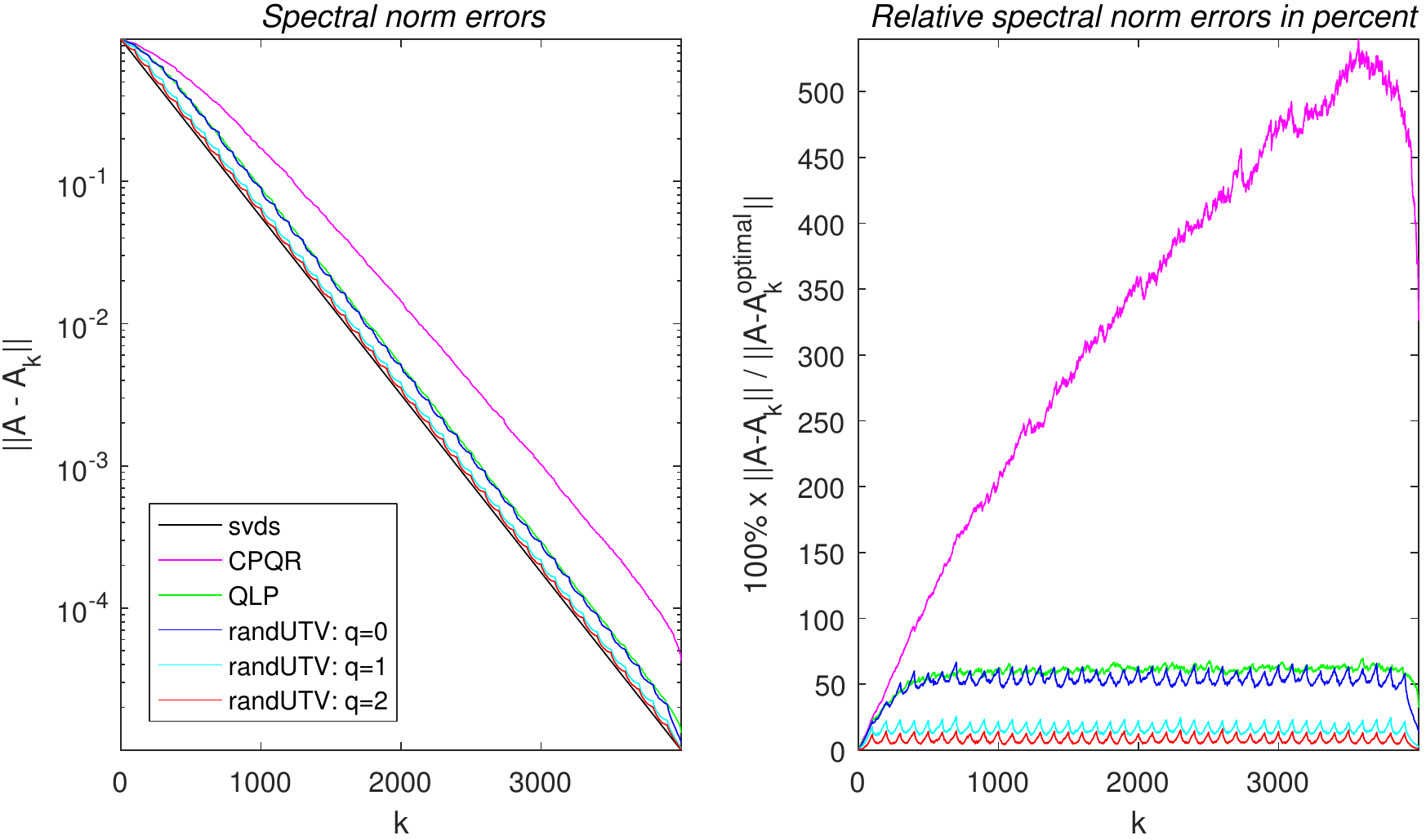}\\
\vspace*{-0.3cm}
\caption{Rank-$k$ approximation errors for the matrix ``Fast Decay''
(see Section \ref{sec:errors}) of size $4000\times 4000$. The block size was $b=100$.
Left: Absolute errors in spectral norm.
The black line marks the theoretically minimal errors. Right: Relative errors, as defined by
(\ref{eq:rel_error}).}
\label{fig:errors_fast_spec}
\end{figure}

\begin{figure}
\vspace*{0.3cm}
\includegraphics[width=0.95\textwidth]{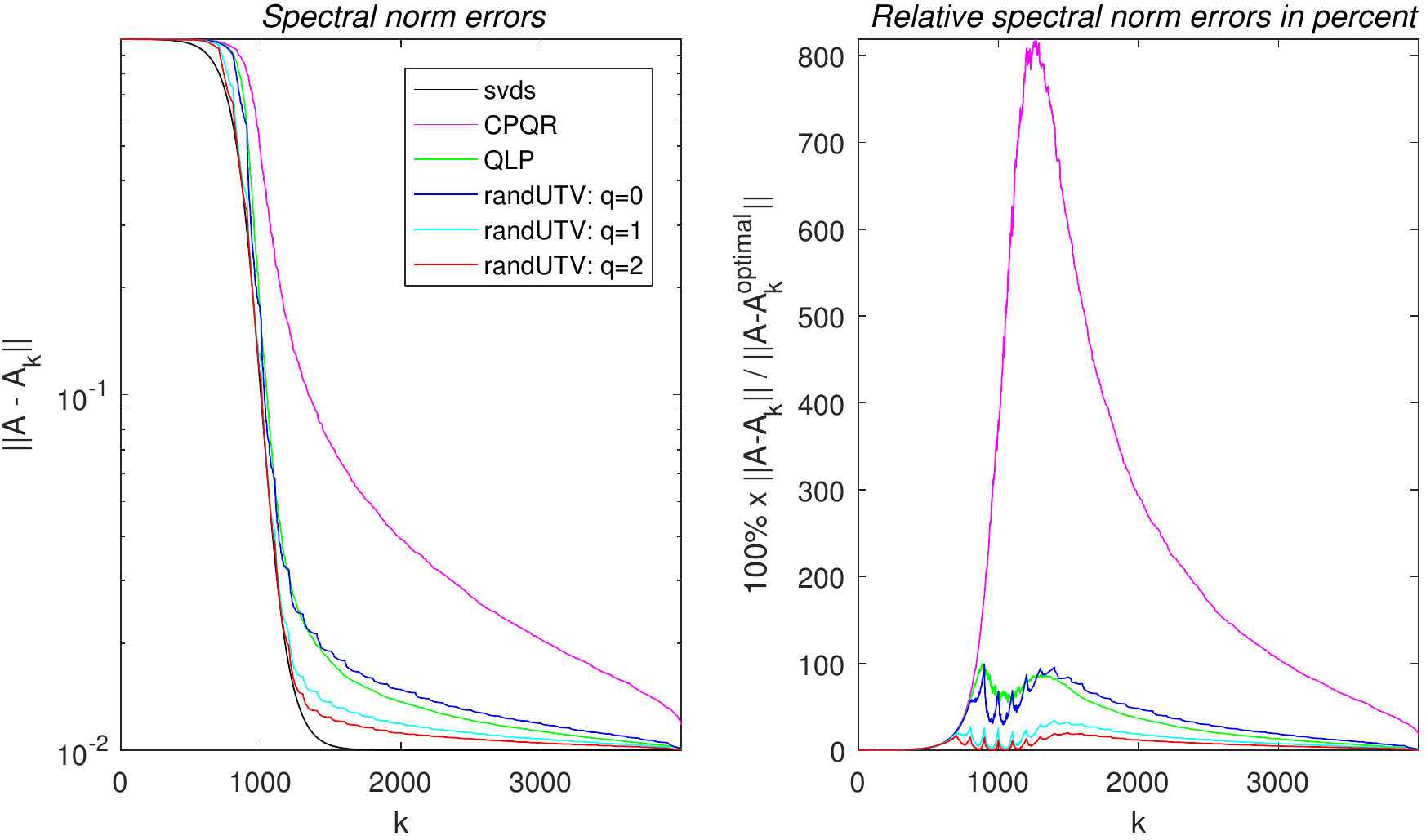}\\
\vspace*{-0.3cm}
\caption{Rank-$k$ approximation errors for the matrix ``S-shaped decay''
(see Section \ref{sec:errors}) of size $4000\times 4000$. The block size was $b=100$.
Left: Absolute errors in spectral norm.
The black line marks the theoretically minimal errors. Right: Relative errors, as defined by
(\ref{eq:rel_error}).}
\label{fig:errors_S_spec}
\end{figure}

\begin{figure}
\vspace*{0.3cm}
\includegraphics[width=0.95\textwidth]{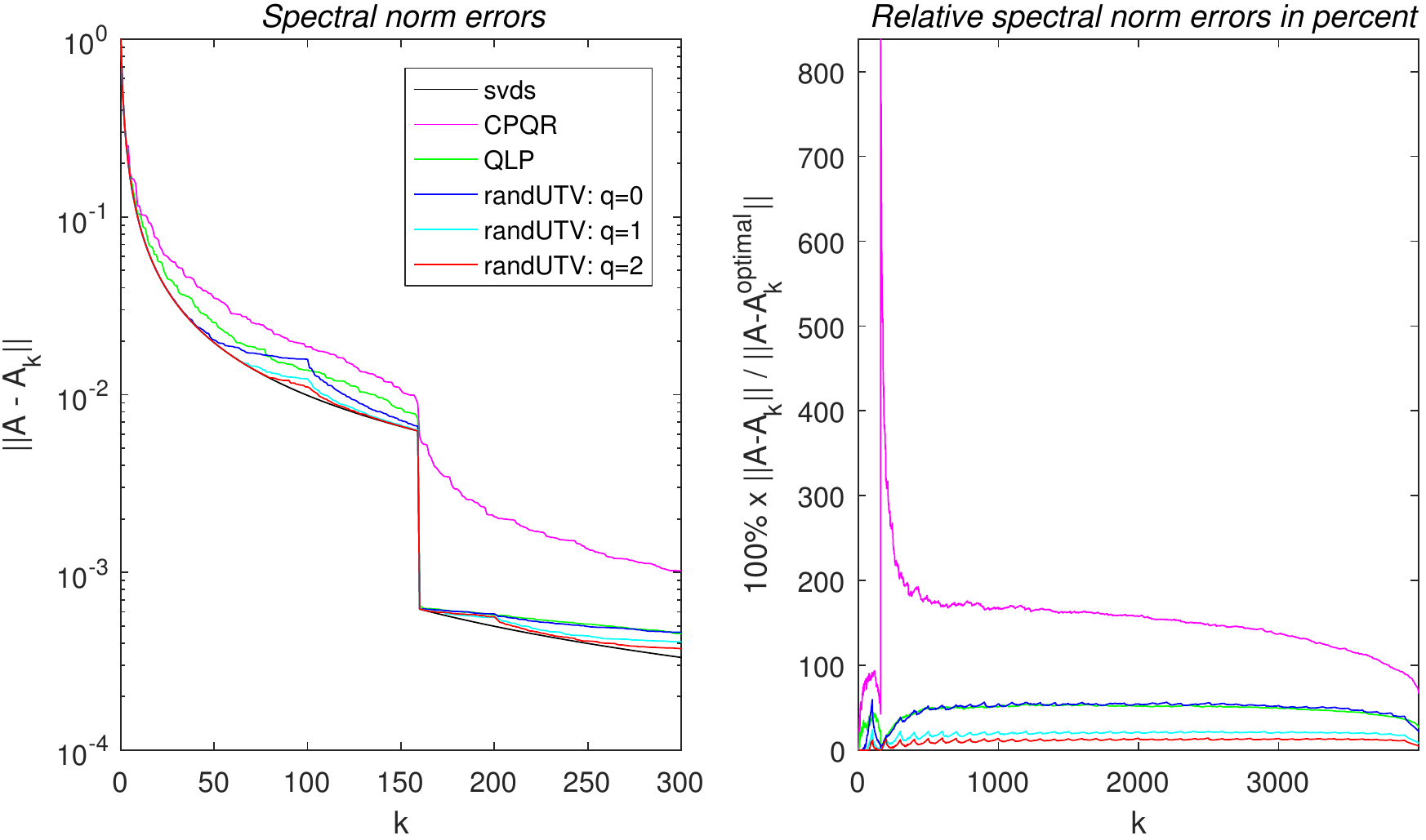}\\
\vspace*{-0.3cm}
\caption{Rank-$k$ approximation errors for the matrix ``Gap''
(see Section \ref{sec:errors}) of size $4000\times 4000$. The block size was $b=100$.
Left: Absolute errors in spectral norm. The black line marks the theoretically minimal errors.
(Observe that we zoomed in on the area around the gap.)
Right: Relative errors, as defined by
(\ref{eq:rel_error}).}
\label{fig:errors_gap_spec}
\end{figure}

\begin{figure}
\vspace*{0.3cm}
\includegraphics[width=0.95\textwidth]{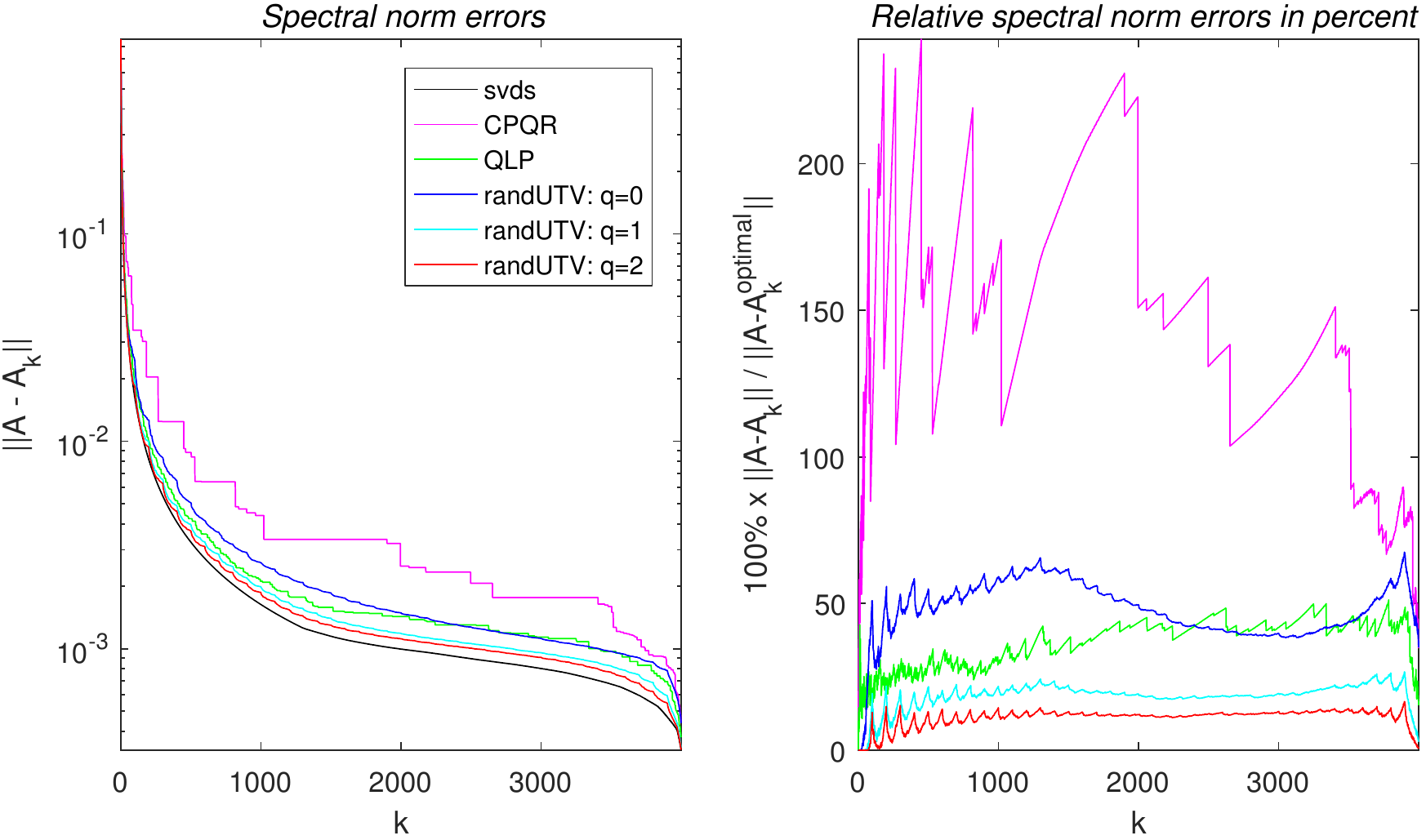}\\
\vspace*{-0.3cm}
\caption{Rank-$k$ approximation errors for the matrix ``BIE''
(see Section \ref{sec:errors}) of size $4000\times 4000$. The block size was $b=100$.
Left: Absolute errors in spectral norm.
The black line marks the theoretically minimal errors. Right: Relative errors, as defined by
(\ref{eq:rel_error}).}
\label{fig:errors_BIE_spec}
\end{figure}

For each test matrix, we computed the error
\begin{equation}
\label{eq:def_error}
e_{k} = \|\mtx{A} - \mtx{A}_{k}\|
\end{equation}
where $\mtx{A}_{k}$ is the rank-$k$ approximation resulting from either of the
three techniques discussed in this manuscript:
\begin{align}
\label{eq:Ak_SVD}
\textit{SVD:}&&     \mtx{A}_{k} =&\ \mtx{A}_{k}^{\rm optimal} = \mtx{U}(:,1:k)\,\mtx{D}(1:k,1:k)\,\mtx{V}(:,1:k)^{*},\\
\label{eq:Ak_CPQR}
\textit{CPQR:}&&    \mtx{A}_{k} =&\ \mtx{Q}(:,1:k)\,\mtx{R}(1:k,:)  \,\mtx{P}^{*},\\
\label{eq:Ak_randUTV}
\textit{randUTV:}&& \mtx{A}_{k} =&\ \mtx{U}(:,1:k)\,\mtx{T}(1:k,:)  \,\mtx{V}^{*}.
\end{align}
For \randUTV{}, we ran the experiment with zero, one, and two steps of power iteration ($q=0,1,2$).
In addition to the direct errors defined by (\ref{eq:def_error}), we also calculated the relative
errors, as defined via
\begin{equation}
\label{eq:rel_error}
e_{k}^{\rm relative} = 100\% \times
\frac{\|\mtx{A} - \mtx{A}_{k}\|}{\|\mtx{A} - \mtx{A}_{k}^{\rm optimal}\|}.
\end{equation}
The results are shown in Figures
\ref{fig:errors_fast_spec},
\ref{fig:errors_S_spec},
\ref{fig:errors_gap_spec}, and
\ref{fig:errors_BIE_spec}.
The figures also report the errors resulting from a UTV factorization
that G.W.~Stewart proposed in \cite{1999_stewart_QLP},
precisely for purposes of low-rank approximation and estimation
of singular values.
To be precise,
Stewart builds a factorization $\mtx{A} = \mtx{U}\mtx{L}\mtx{V}^{*}$ where $\mtx{U}$ and $\mtx{V}$
are orthonormal, and $\mtx{L}$ is lower triangular. The procedure is to first
compute a CPQR factorization of $\mtx{A}$ so that $\mtx{A} = \mtx{Q}_{1}\mtx{R}_{1}\mtx{P}_{1}^{*}$.
Then compute a CPQR of the transpose of the upper triangular matrix $\mtx{R}_{1}$ so that
$\mtx{R}_{1}^{*} = \mtx{Q}_{2}\mtx{L}^{*}\mtx{P}_{2}^{*}$. Finally,  set $\mtx{U} = \mtx{Q}_{1}\mtx{P}_{2}$
and $\mtx{V} = \mtx{P}_{1}\mtx{Q}_{2}$. The rank-$k$ approximant is then
\begin{align}
\label{eq:Ak_QLP}
\textit{QLP:}&&     \mtx{A}_{k} =&\ \mtx{U}\,\mtx{L}(:,1:k)\,\mtx{V}(:,1:k)^{*}.
\end{align}

Based on the errors shown in Figures \ref{fig:errors_fast_spec}--\ref{fig:errors_BIE_spec},
we make several empirical observations:
(1) randUTV is much better than CPQR at computing low-rank approximations.
Even when no power iteration ($q=0$) is used, errors from randUTV are substantially smaller. When one or two steps
of power iteration are taken ($q=1$ or $q=2$), the errors become close to optimal in all cases studied.
(2) For the matrix with a gap in its singular values (cf.~Figure \ref{fig:errors_gap_spec}), randUTV
performs remarkably well in that both $\sigma_{150}$ and $\sigma_{151}$ are approximated to high accuracy.
(3) The relative errors resulting from \randUTV{} are consistently small, and much more reliably small
than those resulting from CPQR.
(4) Comparing the ``QLP'' factorization of Stewart \cite{1999_stewart_QLP} to \randUTV{}, we see that
Stewart's algorithm results in errors that are similar to those resulting from \randUTV{} with $q=0$.
As soon as the power parameter is increased, \randUTV{} tends to perform better. (We observe that in
terms of speed, Stewart's QLP algorithm relies on two CPQR factorizations, which makes it much slower 
than \randUTV{}.)

All error results shown in this section refer to errors measured in the spectral norm. When errors are
measured in the Frobenius norm, \randUTV{} performs even better, as shown in Figures
\ref{fig:errors_fast_frob_detailed},
\ref{fig:errors_S_frob_detailed},
\ref{fig:errors_gap_frob_detailed}, and
\ref{fig:errors_BIE_frob_detailed}.
The effects of including over-sampling in the algorithm, as described in Remark \ref{remark:positivep},
is illustrated in numerical experiments given in Appendix \ref{app:oversampling}. These experiments show
that over-sampling does improve the error, but also that the improvement is almost imperceptible.
For most applications, over-sampling is in our experience not worth the additional effort.

\subsection{Concentration of mass to the diagonal}

\begin{figure}
\includegraphics[width=0.95\textwidth]{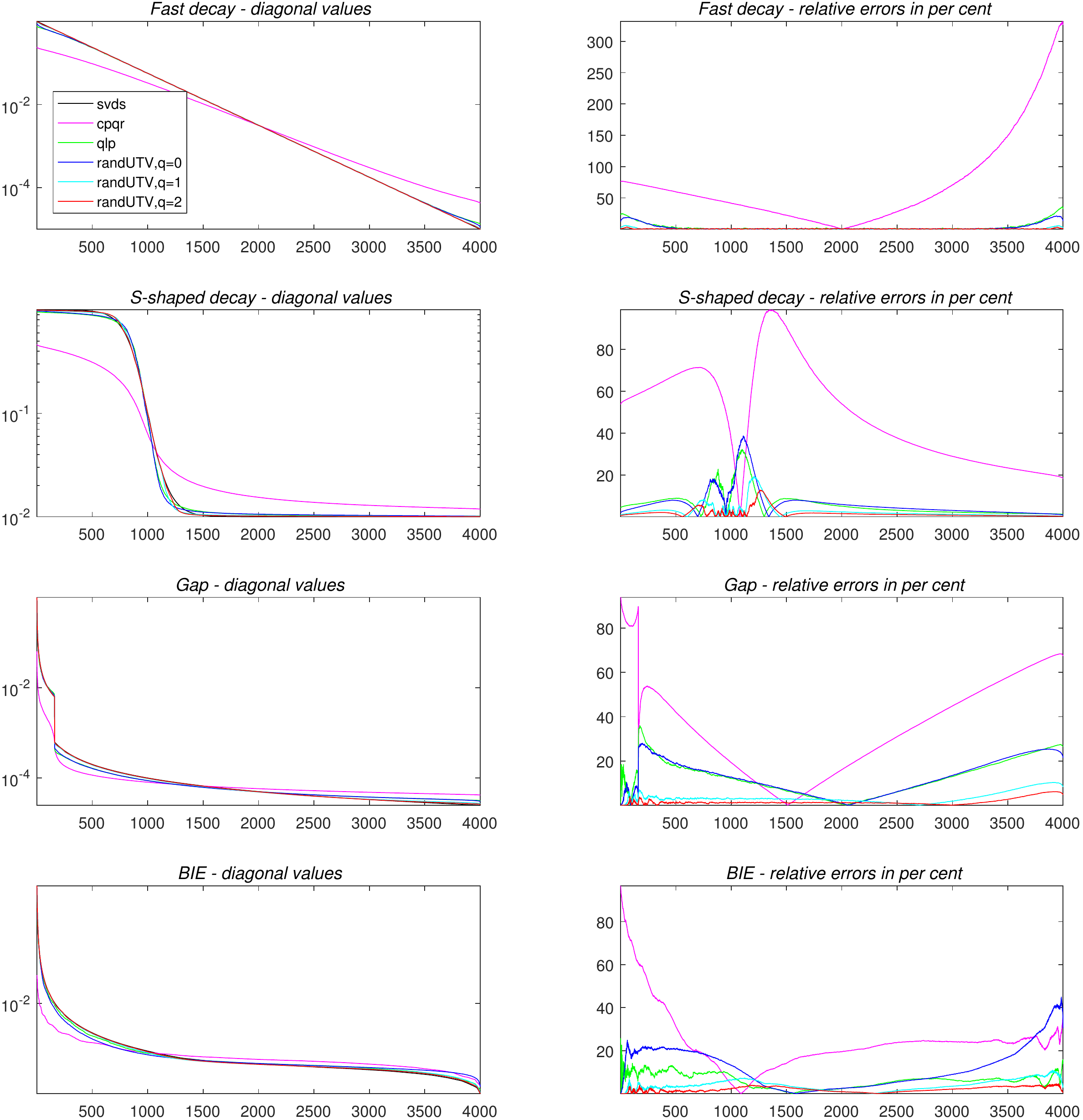}

\caption{The results shown here illustrate how well \randUTV{} approximates the singular values
of four different matrices of size $4000 \times 4000$. For reference, we also include the diagonal
values of the middle factors in a plain CPQR, and in Stewart's ``QLP factorization''  \cite{1999_stewart_QLP}.
The left column shows the diagonal entries themselves. The right column shows relative errors, so that, e.g.,
the magenta line for CPQR plots $100\% \times |\mtx{R}(i,i) - \sigma_{i}|/\sigma_{i}$ versus $i$.}
\label{fig:diags}
\end{figure}

In this section, we investigate our claim that \randUTV{} produces a matrix $\mtx{T}$ whose
diagonal entries are close approximations to the singular values of the given matrix.
(We recall that in the factorization $\mtx{A} = \mtx{U}\mtx{T}\mtx{V}^{*}$, the matrix $\mtx{T}$
is upper triangular, but with the entries above the diagonal very small in magnitude which
forces the diagonal entries to approach the corresponding singular values.)
Figure \ref{fig:diags} shows the results for the four different test matrices
described in Section \ref{sec:errors}, again for matrices of size $4\,000 \times 4\,000$.
For reference, we also show the diagonal entries of the ``R-factor'' in a CPQR, and the diagonal
entries of the ``L-factor'' in Stewart's QLP factorization \cite{1999_stewart_QLP}.
We see that both Stewart's and our algorithm results in far better results than plain
CPQR. \randUTV{} roughly matches the accuracy of the QLP when $q=0$, and does better as
$q$ is increased to one or two. (We recall that \randUTV{} is much faster than the QLP
method.)

\section{Availability of code}
\label{sec:code}

Implementations of the discussed algorithm are available
under 3-clause (modified) BSD license from:

\begin{center}
\texttt{https://github.com/flame/randutv}
\end{center}

\noindent
This repository includes two different implementations:
one to be used with the LAPACK library~\cite{LAPACK3},
and the other one to be used
with the libFLAME~\cite{inverse-siam,FLAME} library.

\textit{\textbf{Acknowledgements:}}
The research reported was supported by DARPA,
under the contract N66001-13-1-4050, and by the NSF,
under the awards DMS-1407340 and DMS-1620472.

\nocite{golub}
\bibliography{main_bib}
\bibliographystyle{amsplain}

\clearpage
\begin{appendix}

\section{Supplementary numerical results}

In this appendix, we present some additional numerical data
that illustrates the performance and accuracy of the algorithms in more detail.

\subsection{Additional timing data}
\label{app:moretimes}

In Section \ref{sec:speed} we provided graphs showing the computation speed of
various factorization algorithms for the case
when all orthonormal matrices are build explicitly.
In this section, we additionally provide the results for the times
required when the orthonormal matrices are not built.
Table~\ref{table:times_noort} 
shows these times when executed on 1, 4, and 14 cores, respectively.
These numbers
were illustrated graphically in the left column of Figure \ref{fig:cputimes}.

\begin{table}
\begin{tabular}{r|r|r|rrr}
\multicolumn{6}{c}{Computational times when executed on a single core} \\
\multicolumn{1}{c|}{$n$} &
\multicolumn{1}{c|}{$T_{\rm svd}$} &
\multicolumn{1}{c|}{$T_{\rm cpqr}$} &
\multicolumn{3}{c}{$T_{\rm randUTV}$} \\
\multicolumn{1}{c|}{} &
\multicolumn{1}{c|}{} &
\multicolumn{1}{c|}{} &
\multicolumn{1}{|c}{$q=0$} &
\multicolumn{1}{c}{$q=1$} &
\multicolumn{1}{c}{$q=2$} \\ \hline
    500 &  3.84e-02 &  1.83e-02 &  4.17e-02 &  4.78e-02 &  5.43e-02 \\
   1000 &  2.57e-01 &  1.24e-01 &  2.02e-01 &  2.48e-01 &  2.93e-01 \\
   2000 &  1.90e+00 &  9.25e-01 &  1.19e+00 &  1.53e+00 &  1.87e+00 \\
   3000 &  3.44e+00 &  4.91e+00 &  3.90e+00 &  5.12e+00 &  6.33e+00 \\
   4000 &  8.48e+00 &  1.35e+01 &  9.13e+00 &  1.21e+01 &  1.51e+01 \\
   5000 &  5.67e+01 &  2.73e+01 &  1.72e+01 &  2.32e+01 &  2.91e+01 \\
   6000 &  1.01e+02 &  4.81e+01 &  2.95e+01 &  3.98e+01 &  5.01e+01 \\
   8000 &  2.50e+02 &  1.20e+02 &  6.95e+01 &  9.37e+01 &  1.18e+02 \\
  10000 &  4.83e+02 &  2.30e+02 &  1.31e+02 &  1.78e+02 &  2.26e+02 \\
\end{tabular}
\vspace*{0.5cm}

\begin{tabular}{r|r|r|rrr}
\multicolumn{6}{c}{Computational times when executed on 4 cores} \\
\multicolumn{1}{c|}{$n$} &
\multicolumn{1}{c|}{$T_{\rm svd}$} &
\multicolumn{1}{c|}{$T_{\rm cpqr}$} &
\multicolumn{3}{c}{$T_{\rm randUTV}$} \\
\multicolumn{1}{c|}{} &
\multicolumn{1}{c|}{} &
\multicolumn{1}{c|}{} &
\multicolumn{1}{|c}{$q=0$} &
\multicolumn{1}{c}{$q=1$} &
\multicolumn{1}{c}{$q=2$} \\ \hline
    500 &  2.60e-02 &  1.31e-02 &  6.54e-02 &  6.85e-02 &  7.24e-02 \\
   1000 &  1.70e-01 &  6.65e-02 &  2.28e-01 &  2.47e-01 &  2.62e-01 \\
   2000 &  9.28e-01 &  4.05e-01 &  8.76e-01 &  1.01e+00 &  1.12e+00 \\
   3000 &  1.72e+00 &  1.21e+00 &  2.06e+00 &  2.47e+00 &  2.81e+00 \\
   4000 &  3.65e+00 &  4.46e+00 &  4.23e+00 &  5.21e+00 &  5.99e+00 \\
   5000 &  1.85e+01 &  1.10e+01 &  7.22e+00 &  9.02e+00 &  1.07e+01 \\
   6000 &  3.49e+01 &  2.03e+01 &  1.16e+01 &  1.46e+01 &  1.76e+01 \\
   8000 &  8.94e+01 &  4.92e+01 &  2.43e+01 &  3.14e+01 &  3.82e+01 \\
  10000 &  1.85e+02 &  9.78e+01 &  4.52e+01 &  5.89e+01 &  7.16e+01 \\
\end{tabular}
\vspace*{0.5cm}

\begin{tabular}{r|r|r|rrr}
\multicolumn{6}{c}{Computational times when executed on 14 cores} \\
\multicolumn{1}{c|}{$n$} &
\multicolumn{1}{c|}{$T_{\rm svd}$} &
\multicolumn{1}{c|}{$T_{\rm cpqr}$} &
\multicolumn{3}{c}{$T_{\rm randUTV}$} \\
\multicolumn{1}{c|}{} &
\multicolumn{1}{c|}{} &
\multicolumn{1}{c|}{} &
\multicolumn{1}{|c}{$q=0$} &
\multicolumn{1}{c}{$q=1$} &
\multicolumn{1}{c}{$q=2$} \\ \hline
    500 &  2.44e-02 &  1.27e-02 &  8.20e-02 &  8.44e-02 &  8.58e-02 \\
   1000 &  1.02e-01 &  5.00e-02 &  2.56e-01 &  2.66e-01 &  2.75e-01 \\
   2000 &  4.14e-01 &  2.16e-01 &  7.98e-01 &  8.56e-01 &  9.00e-01 \\
   3000 &  8.77e-01 &  5.98e-01 &  1.72e+00 &  1.88e+00 &  2.02e+00 \\
   4000 &  1.87e+00 &  3.21e+00 &  3.07e+00 &  3.41e+00 &  3.74e+00 \\
   5000 &  1.36e+01 &  8.74e+00 &  5.00e+00 &  5.64e+00 &  6.28e+00 \\
   6000 &  2.73e+01 &  1.67e+01 &  7.63e+00 &  8.74e+00 &  9.81e+00 \\
   8000 &  7.40e+01 &  4.18e+01 &  1.50e+01 &  1.76e+01 &  2.03e+01 \\
  10000 &  1.55e+02 &  8.33e+01 &  2.58e+01 &  3.10e+01 &  3.62e+01 \\
\end{tabular}
\vspace*{0.2cm}
\caption{Computational times of different factorizations
when executed on one core (top), 4 cores (middle), and 14 cores (bottom).
No orthonormal matrices are built.}
\label{table:times_noort}
\end{table}


\subsection{Improving the accuracy via over-sampling in the randomization step}
\label{app:oversampling}

We recall that the randomized range finder that forms the key novelty of the algorithm
randUTV has been thoroughly investigated and is well understood both empirically and
via rigorous theory, cf.~Section \ref{sec:randrange}. In most of this prior work, it is
standard to perform a small amount of over-sampling, as described in Remark \ref{remark:oversampling}.
In other words, if one seeks a basis that approximately captures the space spanned by
the dominant $k$ singular vectors of a matrix, we draw $k+p$ samples from the range, where
$p$ is a small integer (common choices are $p=5$ or $p=10$). For a single-step method,
such as the simple ``Randomized SVD'' algorithm described in Remark \ref{remark:RSVD},
such over-sampling is essential to get accurate results. For randUTV, it is a simple
matter to incorporate over-sampling in an analogous way. At each step of the algorithm,
we simply draw $p$ extra samples from the range to form an ``extended'' sampling matrix
$\mtx{Y}$ with $k+p$ columns. We then compute the dominant $k$ left singular vectors of
$\mtx{Y}$, and use these as the basis for the column space of the remainder matrix. The
resulting algorithm is given in Figure \ref{fig:stepUTVoversampling}.

Figures \ref{fig:errors_fast_spec_detailed}--\ref{fig:errors_BIE_frob_detailed}
show the errors resulting from randUTV without over-sampling (solid lines) and
with (dotted lines).
The test matrices and the notation are described in Section \ref{sec:errors}.
We see that over-sampling
does provide some benefit in terms of accuracy, as evidenced by the fact that each dotted line
is lower than the corresponding solid line. However, we also see that the difference is very
minor. Overall, it is our opinion that the additional accuracy obtained by over-sampling is
in the present context not worth the additional work. (Avoiding the introduction of an
additional tuning parameter is of course also helpful.)

\begin{figure}[b]
\centering
\begin{minipage}{120mm}
\begin{verbatim}
function [U,T,V] = stepUTV(A,b,q,p)
  G = randn(size(A,1),b+p);
  Y = A'*G;
  for i = 1:q
    Y = A'*(A*Y);
  end;
  [Z,~,~]  = svd(Y,'econ');
  [V,~]    = qr(Z(:,1:b));
  [U,D,W]  = svd(A*V(:,1:b));
  T        = [D,U'*A*V(:,(b+1):end)];
  V(:,1:b) = V(:,1:b)*W;
return
\end{verbatim}
\end{minipage}

\caption{This code illustrates how the accuracy in the single-step UTV factorization can
be improved by incorporating over-sampling, as discussed in Remark \ref{remark:positivep}
and Section \ref{app:oversampling}. The function shown here has exactly the same input
and output parameters as the function \texttt{stepUTV} in Figure \ref{fig:randUTVmatlab},
with the exception that there is one additional input parameter $p$ that specifies how
much oversampling should be done.}
\label{fig:stepUTVoversampling}
\end{figure}


\begin{figure}
\vspace*{0.3cm}
\includegraphics[width=0.95\textwidth]{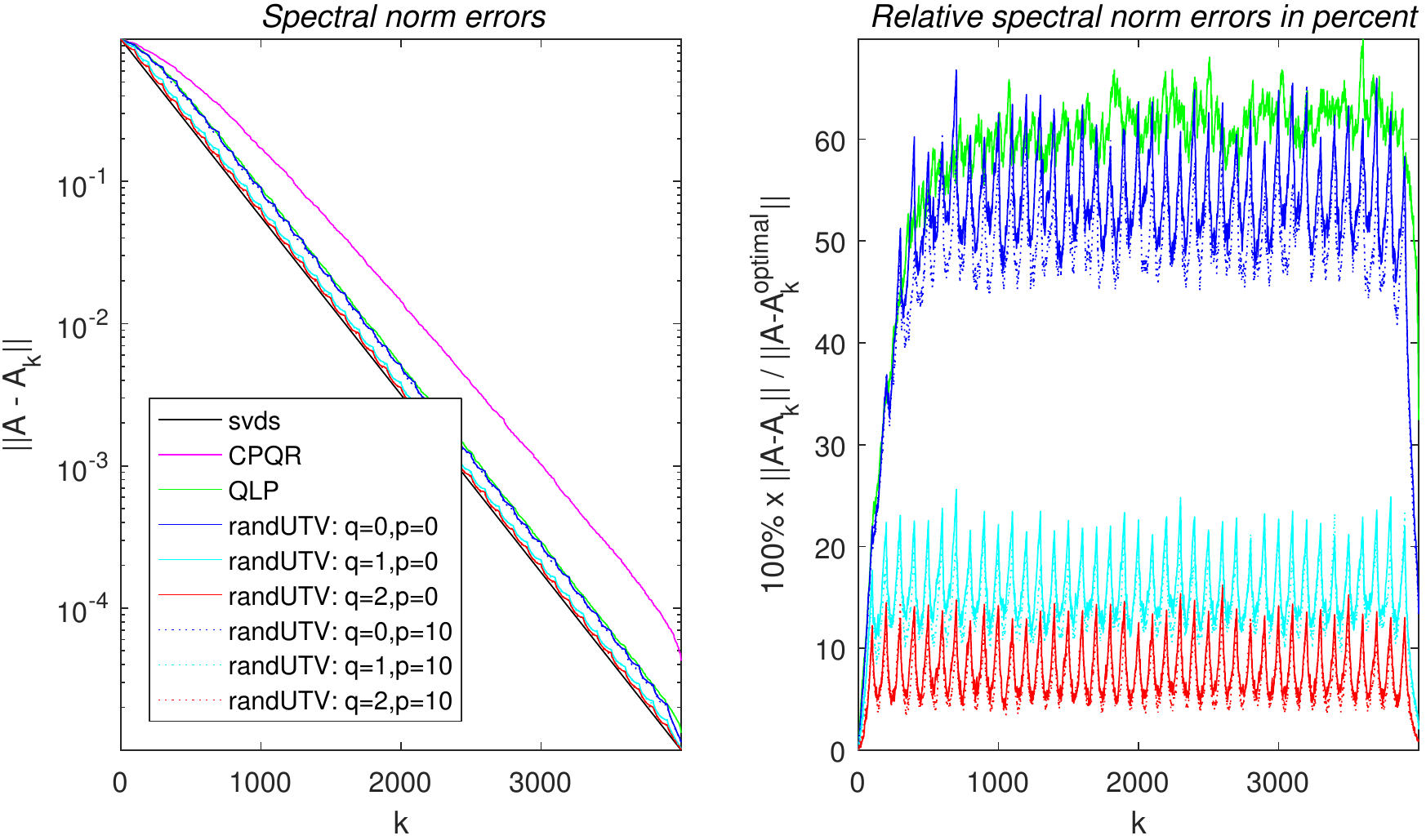}\\
\vspace*{-0.3cm}
\caption{Spectral norm errors for the ``Fast Decay'' matrix.
This figure is identical to Figure \ref{fig:errors_fast_spec}, except that we
now include lines showing the effect of over-sampling (dotted lines).}
\label{fig:errors_fast_spec_detailed}
\end{figure}

\begin{figure}
\vspace*{0.3cm}
\includegraphics[width=0.95\textwidth]{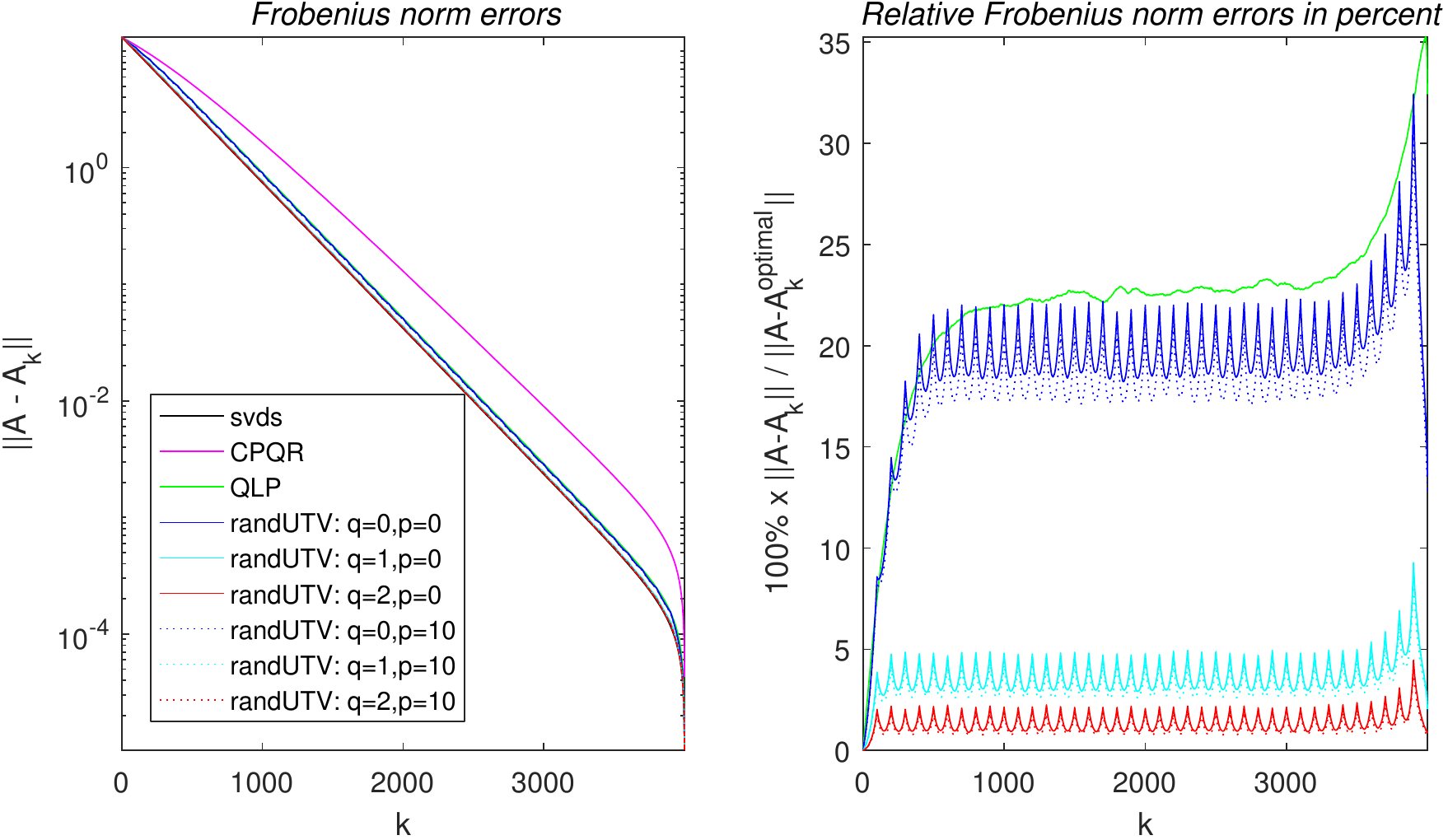}\\
\vspace*{-0.3cm}
\caption{Frobenius norm errors for the ``Fast Decay'' matrix.
This figure is identical to Figure \ref{fig:errors_fast_spec_detailed},
except that errors are now measured in the Frobenius norm.}
\label{fig:errors_fast_frob_detailed}
\end{figure}


\begin{figure}
\vspace*{0.3cm}
\includegraphics[width=0.95\textwidth]{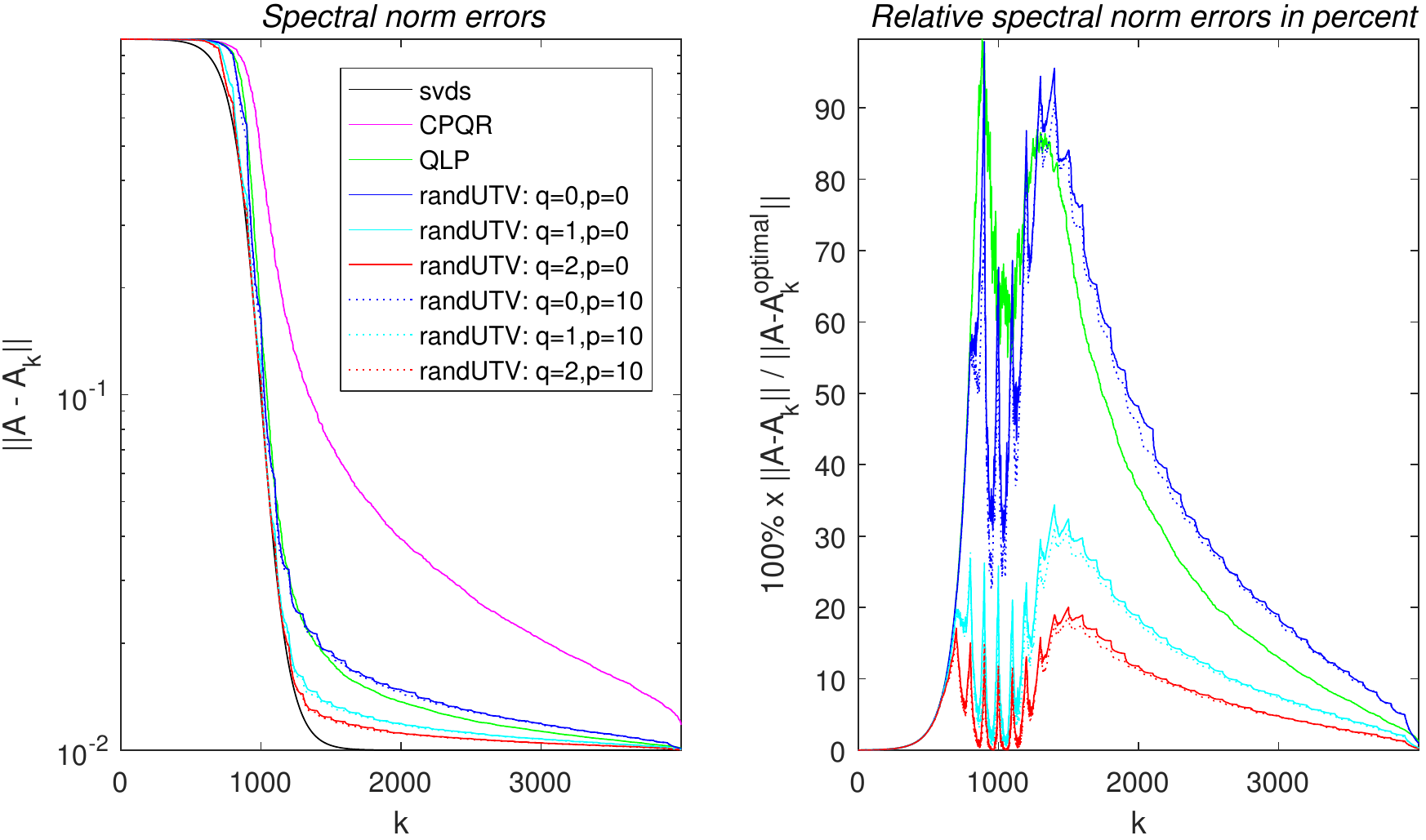}\\
\vspace*{-0.3cm}
\caption{Spectral norm errors for the ``S-shaped Decay'' matrix.
This figure is identical to Figure \ref{fig:errors_S_spec}, except that we
now include lines showing the effect of over-sampling (dotted lines).}
\label{fig:errors_S_spec_detailed}
\end{figure}

\begin{figure}
\vspace*{0.3cm}
\includegraphics[width=0.95\textwidth]{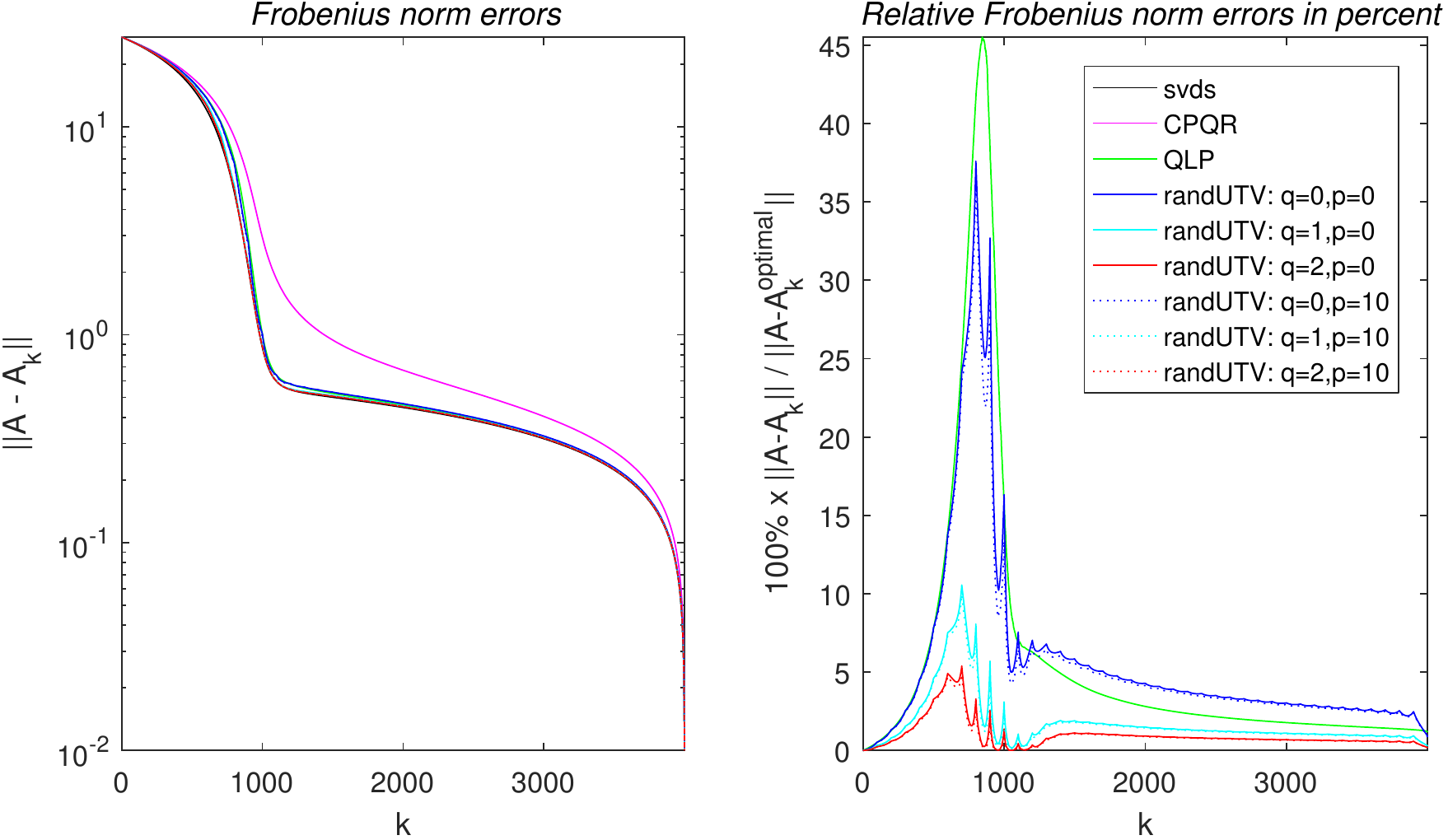}\\
\vspace*{-0.3cm}
\caption{Frobenius norm errors for the ``S-shaped Decay'' matrix.
This figure is identical to Figure \ref{fig:errors_S_spec_detailed},
except that errors are now measured in the Frobenius norm.}
\label{fig:errors_S_frob_detailed}
\end{figure}


\begin{figure}
\vspace*{0.3cm}
\includegraphics[width=0.95\textwidth]{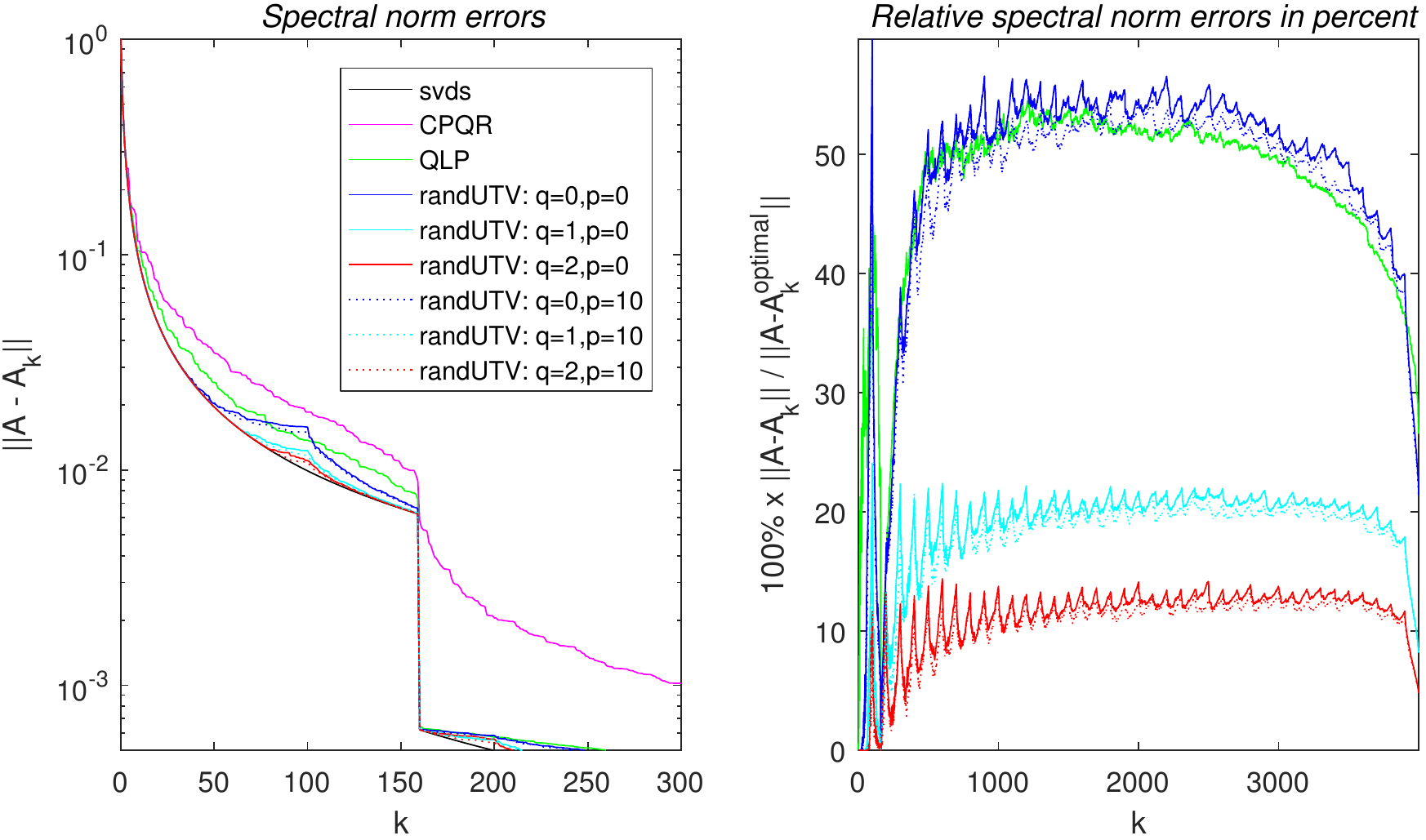}\\
\vspace*{-0.3cm}
\caption{Spectral norm errors for the ``Gap'' matrix.
This figure is identical to Figure \ref{fig:errors_gap_spec}, except that we
now include lines showing the effect of over-sampling (dotted lines).}
\label{fig:errors_gap_spec_detailed}
\end{figure}

\begin{figure}
\vspace*{0.3cm}
\includegraphics[width=0.95\textwidth]{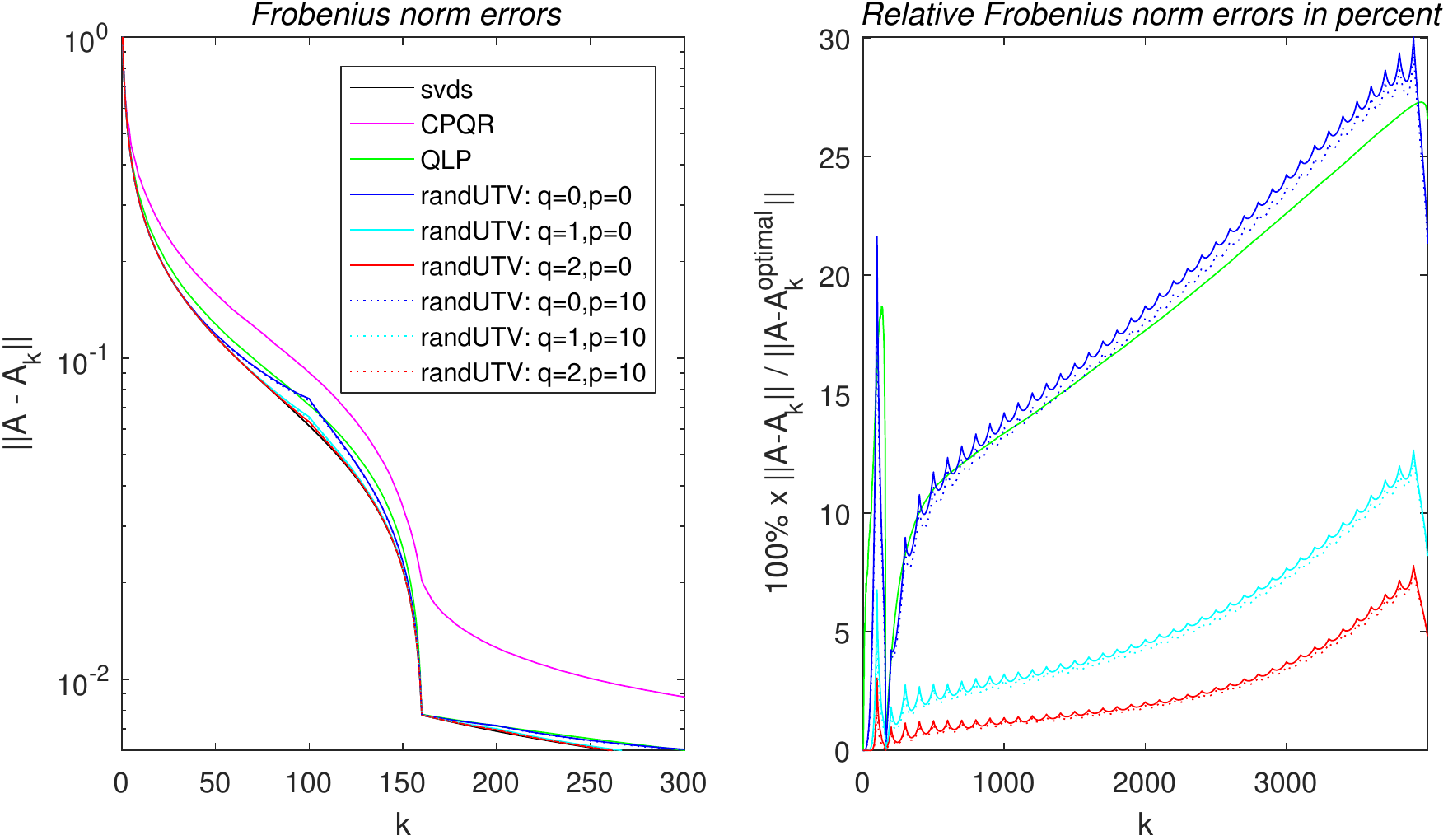}\\
\vspace*{-0.3cm}
\caption{Frobenius norm errors for the ``Gap'' matrix.
This figure is identical to Figure \ref{fig:errors_gap_spec_detailed},
except that errors are now measured in the Frobenius norm.}
\label{fig:errors_gap_frob_detailed}
\end{figure}


\begin{figure}
\vspace*{0.3cm}
\includegraphics[width=0.95\textwidth]{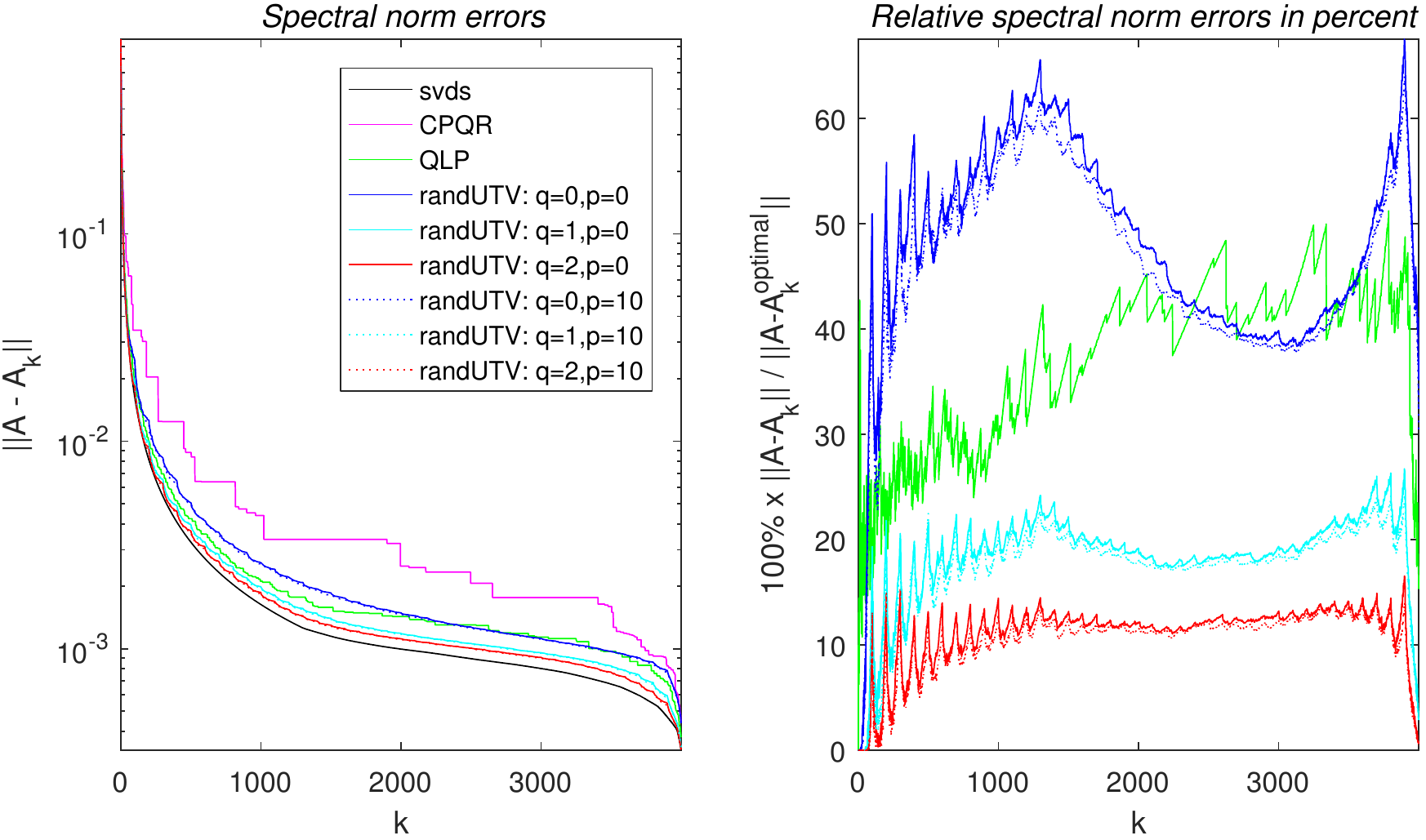}\\
\vspace*{-0.3cm}
\caption{Spectral norm errors for the ``BIE'' matrix.
This figure is identical to Figure \ref{fig:errors_BIE_spec}, except that we
now include lines showing the effect of over-sampling (dotted lines).}
\label{fig:errors_BIE_spec_detailed}
\end{figure}

\begin{figure}
\vspace*{0.3cm}
\includegraphics[width=0.95\textwidth]{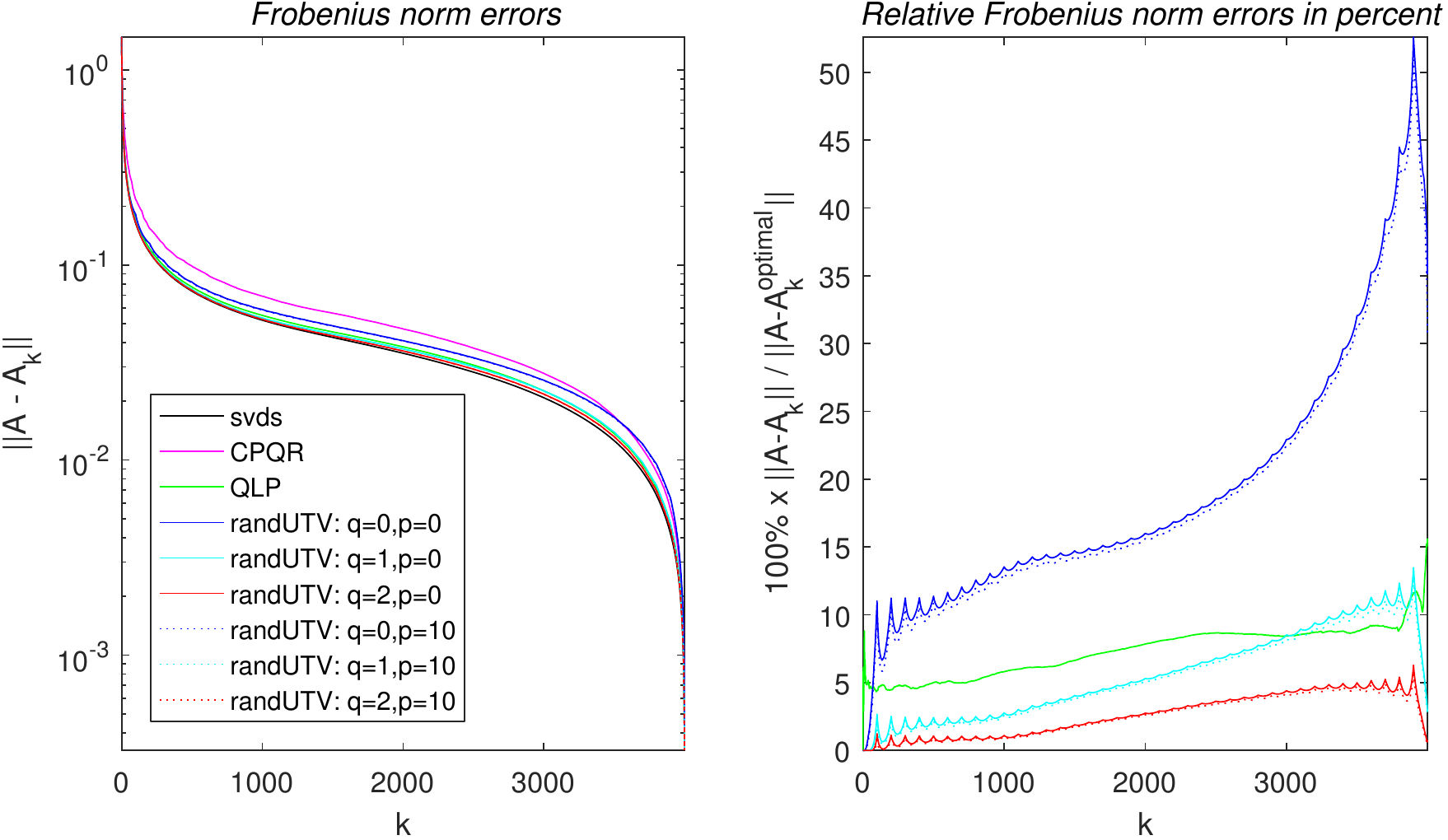}\\
\vspace*{-0.3cm}
\caption{Frobenius norm errors for the ``BIE'' matrix.
This figure is identical to Figure \ref{fig:errors_BIE_spec_detailed},
except that errors are now measured in the Frobenius norm.}
\label{fig:errors_BIE_frob_detailed}
\end{figure}

\end{appendix}

\end{document}